\documentclass[12pt,a4paper,draft]{article}

\usepackage[reqno]{amsmath}
\usepackage{amssymb,amsthm,upref,amscd}
\usepackage[T1]{fontenc}
\usepackage{times}

\newtheorem{theorem}{Theorem}[section]

\newtheorem{definition}[theorem]{Definition}

\newtheorem{lemma}[theorem]{Lemma}
\newtheorem{proposition}[theorem]{Proposition}
\newtheorem{remark}[theorem]{Remark}

\theoremstyle{definition} \theoremstyle{remark}
\numberwithin{equation}{section}
\begin{document}

\title{A proof for the conjecture on superlinear problems with
Ambrosetti-Rabinowitz condition \\
[3mm] {\footnotesize Dedicated to Professor Paul Rabinowitz on the occasion
of his 87th birthday }}
\author{Chong Li \ and \ Shujie Li \thanks{%
The first author is supported by NSFC(11871066), AMSS(E3550105). The second
author is supported by NSFC(11871066). }}
\date{}
\maketitle

\begin{abstract}
This paper is devoted to exploring a new minimax approach by introducing a
characteristic mapping family which is invariant under the smooth descending
flow for initial value. The minimax approach is self-contained, and its
features are markedly different from standard ones, as it identifies the
existence of critical points and intrinsically presents a lower-bound
estimate for the generalized Morse index at the corresponding critical
point. This quantity can be effectively viewed as an alternative to the
group action. As applications, under the Ambrosetti-Rabinowitz condition we
offer a positive answer to the long-standing open problem on the existence
of infinitely many distinct solutions for superlinear elliptic equations
without symmetric hypothesis.
\end{abstract}

\section{Introduction}

Originating in the pioneer work of Ambrosetti-Rabinowitz \cite{AR} in 1973,
followed by Rabinowitz \cite{Ra1}, Benci-Rabinowitz \cite{BR}, Ni \cite{Ni},
Chang \cite{Cha1}, Chang \cite{Cha2}, Wang \cite{W} and others, minimax
methods have been powerful tools for finding saddle points of superquadratic
functionals. In particular, from its first appearance in the work of
Ambrosetti and Rabinowitz \cite{AR} the existence and multiplicity of
solutions to the nonlinear elliptic equation:

\begin{equation}
\left\{
\begin{array}{cc}
-\Delta u=f\left( x,u\right) , & x\in \Omega , \\
u=0, & x\in \partial \Omega ,
\end{array}
\right.  \label{eq1}
\end{equation}
$\Omega $ a bounded domain with smooth boundary $\partial \Omega $ in $\Bbb{R%
}^N$, has been widely studied by many authors under the following
assumptions:

$\left( f_1\right) $ $f\in C\left( \overline{\Omega }\times \mathbb{R},%
\mathbb{R}\right) $, there are $c>0$ and $p\in \left( 2,2^{*}\right) $ such
that $\left| f\left( x,t\right) \right| \leq c\left( 1+\left| t\right|
^{p-1}\right) $ for all $x\in \Omega $ and $t\in \mathbb{R}$, where $2^{*}=%
\frac{2N}{N-2}$ for $N\geq 3$ and $2^{*}=\infty $ for $N=1$,$2$;

$\left( f_2\right) $ (the Ambrosetti-Rabinowitz condition) There are $\theta
>2$ and $M>0$, such that for all $x\in \Omega $ and $\left| t\right| \geq M$%
,
\begin{equation}
0<\theta F\left( x,t\right) \leq tf\left( x,t\right) ,  \label{eq2}
\end{equation}
where $F\left( x,t\right) =\int_0^tf\left( x,s\right) ds$;

$\left( f_3\right) $ $f\left( x,t\right) =o\left( t\right) $ uniformly in $x$
as $t\rightarrow 0$.

$\left( f_4\right) $ $f\left( x,-t\right) =-f\left( x,t\right) $, $\forall
x\in \Omega $, $\forall t\in \Bbb{R}$.

Let $X=H_0^1\left( \Omega \right) $ be equipped with the norm
\[
\left\| u\right\| _X=\left( \int_\Omega \left| \nabla u\right| ^2\right)
^{\frac 12}
\]
and with the inner product
\[
\left\langle u,v\right\rangle _X=\int_\Omega \nabla u\nabla vdx,
\]
then the weak solutions of $\left( \ref{eq1}\right) $ are the critical
points of the following energy functional
\[
J\left( u\right) =\frac 12\int_\Omega \left| \nabla u\right|
^2dx-\int_\Omega F\left( x,u\right) dx.
\]

Based on the assumptions $\left( f_1\right) $--$\left( f_3\right) $
Ambrosetti and Rabinowitz showed that the equation $\left( \ref{eq1}\right) $
admits at least two nontrivial solutions in \cite{AR}. Under the same
assumptions, a major breakthrough was made by Wang \cite{W}, who found the
third one in 1991 by employing linking approach and Morse theory.

Furthermore, under the hypotheses $\left( f_1\right) $--$\left( f_4\right) $%
, Ambrosetti and Rabinowitz \cite{AR} proved the existence of infinitely
many solutions. A fundamental open problem, independently posed by
Rabinowitz and Struwe, is whether the equation $\left( \ref{eq1}\right) $\
still admits infinitely many distinct solutions once the oddness assumption $%
\left( f_4\right) $ is removed. It has been widely expected that the answer
to this problem is positive. Indeed, this is true for the case $N=1$(see
Nehari \cite{Ne}). However, to the best of our knowledge, for the case $%
N\geq 2$, the renowned conjecture has remained open for over sixty years,
since the work of Nehari in 1961(see Struwe \cite[p.124]{Stru3} for
details). The literature on superlinear elliptic equations is enormous and
the readers may consult Struwe \cite{Stru1}, Struwe \cite{Stru2}, Struwe
\cite{Stru3}, Turner \cite{Turner}, Rabinowitz \cite{Ra2}, Rabinowitz \cite
{Ra3}, Rabinowitz \cite{Ra4}, Bartsch-Chang-Wang \cite{BCW}, Bartsch-Weth
\cite{BW}, Rabinowitz-Su-Wang \cite{RSW}, Ekeland-Ghoussoub \cite{EG},
Ghoussoub-Fang \cite{FG}, Szulkin \cite{Szu}, Li-Liu \cite{LLZ}, Li-Li \cite
{LiLi}, Li-Liu \cite{LiLiu}, Sun-Su-Tian \cite{SST} and references therein.

Attempts to resolve this problem over decades have faced significant
challenges, as standard minimax theories are not directly applicable due to
the absence of the symmetry or topological linking structures they crucially
rely on. The current paper investigates the question from a different
perspective and adopts an alternative approach. In order to circumvent the
substantial difficulty posed by standard theories, we develop a new minimax
method by introducing a distinctive mapping family. We show that smoothing
of nonlinear functions induces a smooth descending flow about initial
values, and the family is invariant under the flow. The minimax method
presented here differs from standard ones by incorporating a smooth mapping
family and providing explicit lower bounds for the generalized Morse
indices, and enables us to finally settle the long-standing conjecture on
the existence of infinitely many distinct solutions for superlinear elliptic
equations under the Ambrosetti-Rabinowitz condition alone, without the
oddness hypothesis $\left( f_4\right) $.

Let $E$ be a Hilbert space and $\Phi \in C^2\left( E,\Bbb{R}\right) $. For
the sake of convenience, in what follows, we denote by $M\left( \Phi
,u\right) $ the generalized Morse index of $\Phi $ at $u$, i.e., the sum of
the dimensions of the negative and null eigenspace of the operator $\Phi
^{\prime \prime }\left( u\right) $ acting on $E$ if $u$ is a critical point
of $\Phi $ on $E$.

We need the following assumption:$\ $

$\left( f_1^{\prime }\right) $ $f\in C^1\left( \overline{\Omega }\times %
\mathbb{R},\Bbb{R}\right) $, there are $C>0$ and $p\in \left( 2,2^{*}\right)
$ such that
\begin{equation}
\left| f_t^{\prime }\left( x,t\right) \right| \leq C\left( 1+\left| t\right|
^{p-2}\right) ,  \label{eq3}
\end{equation}
for all $x\in \Omega $ and $t\in \Bbb{R}$, where $2^{*}=\frac{2N}{N-2}$ if $%
N\geq 3$ and $2^{*}=\infty $ if $N=1$,$2$.

It is well known that the assumption $\left( f_1^{\prime }\right) $
guarantees $J\in C^2\left( X,\Bbb{R}\right) $. In this paper we give a
positive answer to the above question. Our first main theorem concerning the
existence of multiple solutions of superlinear elliptic equations reads:

\begin{theorem}
\label{thm1.1} Suppose $\left( f_1^{\prime }\right) \left( f_2\right) \left(
f_3\right) $ hold, then the functional $J$ possesses an unbounded sequence $%
\left\{ u_k\right\} _{k=1}^\infty $ of critical points satisfying

$\left( i\right) $ $M\left( J,u_k\right) \rightarrow \infty $ as $%
k\rightarrow \infty $.

$\left( ii\right) $ $J\left( u_k\right) \rightarrow \infty $ as $%
k\rightarrow \infty $.
\end{theorem}

Let
\[
J_{p,q}\left( u\right) :=\frac 12\int_\Omega \left| \nabla u\right|
^2dx-\frac 1p\int_\Omega \left| u\right| ^pdx-\frac 1q\int_\Omega \left|
u\right| ^{q-1}udx.
\]
As an application of Theorem \ref{thm1.1}, we give a positive answer to the
following pending issue:

\begin{theorem}
\label{thm1.2} The equation
\begin{equation}
\left\{
\begin{array}{cc}
-\Delta u=\left| u\right| ^{p-2}u+\left| u\right| ^{q-1}, & x\in \Omega , \\
u=0, & x\in \partial \Omega ,
\end{array}
\right.  \label{eq4}
\end{equation}
possesses an unbounded solution sequence $\left\{ u_k\right\} _{k=1}^\infty $
with $M\left( J_{p,q},u_k\right) \rightarrow \infty $, $J_{p,q}\left(
u_k\right) \rightarrow \infty $ as $k\rightarrow \infty $, where $2\leq q<p<%
\frac{2N}{N-2}$ if $N\geq 3$, and $2\leq q<p<\infty $ if $N=1$,$2$.
\end{theorem}

\begin{remark}
\label{Remark1.3} Our research focuses on proving the existence of
infinitely many solutions for superlinear problems under the
Ambrosetti-Rabinowitz condition, a direction distinct from the
symmetry-breaking perspective. This work settles the long-standing
conjecture for the general equation $-\Delta u=f\left( x,u\right) $ under
assumptions $\left( f_1^{\prime }\right) $-$\left( f_3\right) $. Our setting
differs fundamentally from the extensively studied model with an external
force $-\Delta u=\left| u\right| ^{p-2}u+f\left( x\right) $, where the
nonlinearity has the specific structure of a homogeneous principal part plus
an additive perturbation. Seminal results for this latter model were
obtained by Bahri-Berestycki and Bahri-Lions. For comprehensive accounts of
that theory, we refer to the dedicated chapter `Multiple critical points of
perturbed symmetric functionals' in \cite{Ra4} and the monograph
\cite[Chapter II.7]{Stru3}. Additionally, our work does not address the
existence of infinitely many sign-changing solutions.
\end{remark}

A sketch of the proof of the main theorem is furnished in this section due
to the intricacy of the argument. The proofs of Theorems \ref{thm1.1} and
\ref{thm1.2} are structured around a novel minimax argument that
incorporates the estimates of explicit lower bounds of the generalized Morse
index directly into seeking for critical points of the energy functional.
The proof outline is divided into five parts and proceeds as follows:

$\left\langle i\right\rangle $ A New Minimax Framework. The method we employ
here to substantiate the existence of critical values for functionals with a
mountain pass structure is straightforward, but differs fundamentally from
the standard one for the proof of the mountain pass theorem. Our primary
objective is to distinguish these critical points by estimating their
generalized Morse indices. Standard minimax theories, which usually rely on
a family of mappings with a group action (e.g., G-homotopic, or
G-cohomotopic mappings, see Ghoussoub \cite{Ghou} for details) or on various
forms of topological linking, have been successfully applied to prove the
multiplicity and provide generalized Morse index(or Morse index) information
for solutions of various problems. However, the standard methods are not
suitable for the broader class of superlinear elliptic problems considered
in this paper, which satisfy the Ambrosetti-Rabinowitz condition but lack
symmetry and do not readily admit the verification of nontrivial topology.

$\left\langle ii\right\rangle $ Smoothing and Finite-Dimensional
Approximation. To obtain the desired lower bounds for the generalized Morse
index of critical points of the energy functional $J$ on $X$, we mollify the
nonlinearity $f$ and obtain smoothed functionals $J_m$. For fixed $m$, we
then investigate critical points of $J_m$ restricted to finite-dimensional
subspaces $X_n$, named by $J_{m,n}$.

$\left\langle iii\right\rangle $ The Invariant Minimax Family and Smooth
Flow. The mapping family for $J_m$ is defined on an annular region of $X_k$%
(with $k<n$), denoted by $D_{r_k,R_k}$. This family is invariant under the
descending flow of $X_n$. The key point is to establish the $C^2$ smoothness
of the flow $\tau \left( t,u\right) $ defined on $\Bbb{R}^{+}\cup \left\{
0\right\} \times X_n$, with respect to the initial value $u\in X_n$. In
Lemma \ref{lem2.1}, the Gronwall-Bellman inequality first ensures the
uniqueness of the flow. This, combined with the $C^2$ smoothness of the
governing vector field, then establishes the injectivity of its derivative
and the flow's $C^2$ smoothness in the initial data(see Lemma \ref{lem2.1}
for details). The mapping family is indeed a restriction of the flow on $%
\Bbb{R}^{+}\cup \left\{ 0\right\} \times \overline{D_{r_k,R_k}}$.
Consequently, for any fixed $t>0$, the flow defines a diffeomorphism from $%
D_{r_k,R_k}$ onto a $C^2$ submanifold $A_t$ of $X_n$.

$\left\langle iv\right\rangle $ Critical Points with Index Bounds. We prove
that the minimax value $C_k^{\left( m,n\right) }$ associated with this
family is a critical value of $J_m$ on $X_n$. Furthermore, by analyzing a
sequence of maximum points $v_{t_j}$ on $A_{t_j}$ for $t_j\rightarrow \infty
$, we obtain a Palais-Smale sequence for $J_m$ on $X_n$. The sum of the
dimensions of the negative and zero subspaces of the Hessian $J_m^{\prime
\prime }\left( v_{t_j}\right) $ is at least $k$, which is ensured by the
triviality of the kernel of $\tau _u^{\prime }\left( t,u\right)
:X_k\rightarrow T_{\tau \left( t,u\right) }A_t$(see also Lemma \ref{lem2.1}%
). Here $T_{\tau \left( t,u\right) }A_t$ is the tangent space of $A_t$ at $%
\tau \left( t,u\right) $. This allows us to conclude the existence of a
critical point $u_{m,n,k}$ of $J_m$ on $X_n$ with the generalized Morse
index $M\left( J_{m,n},u_{m,n,k}\right) \geq k$.

$\left\langle v\right\rangle $ Limit Passage.\textbf{\ }A final two-step
limit process, first as $m\rightarrow \infty $ and then as $n\rightarrow
\infty $, employing the (PS)$^{\text{**}}$ (see Definition \ref{def4.3}) and
(PS)$^{\text{*}}$ conditions respectively, yields a sequence of critical
points $\left\{ u_k\right\} _{k=1}^\infty $ for the functional $J$
satisfying $M\left( J,u_k\right) \geq k$, and thus completes the proof of
the main theorems.

The remainder of this manuscript is organized as follows. Section 2 is
covered by the elaboration on the proof of smoothness of flow. Section 3
addresses the smooth approximation of the nonlinear term. The demonstrations
of Theorems \ref{thm1.1} and \ref{thm1.2} are established in Section 4.

\textbf{Acknowledgments.} We are greatly indebted to Professors Thomas
Bartsch, Andrzej Szulkin and Zhiqiang Wang for various valuable suggestions.
The first author is supported by National Natural Science Foundation of
China(Grant No.11871066), and also by Academy of Mathematics and Systems
Science(Grant No.E3550105). Finally, Chong Li is sincerely grateful for his
family's work behind the scenes over the years.

\section{Smoothness of flow for initial values}

Let $X$ be a Banach space and make the following hypothesis:

$(g_1)$ $g:X\rightarrow X$ is a $C^2$ mapping.

Aiming at concluding the proof of Theorem \ref{thm1.1}, we are motivated by
the following lemma, which serves as a cornerstone of the argument.

\begin{lemma}
\label{lem2.1} Assume $\left( g_1\right) $ holds, i.e., $g:X\rightarrow X$
is a $C^2$ mapping. Then the Cauchy problem
\begin{equation}
\left\{
\begin{array}{c}
\frac{d\tau \left( t,u\right) }{dt}=g\left( \tau \left( t,u\right) \right) ,%
\text{ }\left( t,u\right) \in \Bbb{R}^{+}\times X, \\
\tau \left( 0,u\right) =u\in X,
\end{array}
\right.  \label{eq5}
\end{equation}
has a unique solution $\tau \left( t,u\right) \in C^2\left( \Bbb{R}%
^{+}\times X,X\right) $. Furthermore, the derivative with respect to the
initial value, which we denote by $\tau _u^{\prime }\left( t,u\right) $, as
a bounded linear operator from $X$ to $X$, is injective for all $t$, i.e., $%
\ker \tau _u^{\prime }\left( t,u\right) =\left\{ 0\right\} $. Finally, for
each initial point $u\in X$, if $g\left( \tau \left( t,u\right) \right) =0$
and $g\left( u\right) \neq 0$, then $t=\infty $.
\end{lemma}

For convenience, we briefly review the Gronwall-Bellman inequality, adapted
to our needs:

\begin{proposition}
\label{prop2.2} Suppose that the continuous functions $z\left( t\right) $, $%
b\left( t\right) \geq 0$ for $t\in \left[ t_0,t_1\right] $, and $a\geq 0$
satisfy
\begin{equation}
z\left( t\right) \leq a+\int_{t_0}^tb\left( s\right) z\left( s\right) ds,
\label{eq6}
\end{equation}
then
\begin{equation}
z\left( t\right) \leq a\cdot \exp \left\{ \int_{t_0}^{t_1}b\left( s\right)
ds\right\}  \label{eq7}
\end{equation}
for $t\in \left[ t_0,t_1\right] $. Moreover, if
\begin{equation}
z\left( t\right) \leq z\left( t_0\right) +\int_{t_0}^tb\left( s\right)
z\left( s\right) ds,  \label{eq8}
\end{equation}
then
\begin{equation}
z\left( t_0\right) \cdot \exp \left\{ -\int_{t_0}^tb\left( s\right)
ds\right\} \leq z\left( t\right) \leq z\left( t_0\right) \cdot \exp \left\{
\int_{t_0}^tb\left( s\right) ds\right\}  \label{eq9}
\end{equation}
for $t\in \left[ t_0,t_1\right] $.
\end{proposition}

\textit{Proof of Lemma 2.1. }The first assertion is standard(see, e.g.,
\cite[Chapter IV.1]{Lang}), but our proof is entirely different. Given two
Banach spaces $X$ and $Y$, we denote by $L\left( X,Y\right) $ the set of all
bounded linear operators from $X$ to $Y$ and shortened by $L(X)$ if $X=Y$.
Notice that $\left( \ref{eq5}\right) $ is equivalent to
\begin{equation}
\tau \left( t,u\right) =u+\int_0^tg\left( \tau \left( s,u\right) \right) ds,
\label{eq10}
\end{equation}
If $\tau _u^{\prime }\left( t,u\right) $ exists, then by differentiating
both sides of equation $\left( \ref{eq10}\right) $ with respect to $u$, we
obtain
\begin{equation}
\tau _u^{\prime }\left( t,u\right) =I+\int_0^tg^{\prime }\left( \tau \left(
s,u\right) \right) \tau _u^{\prime }\left( s,u\right) ds,  \label{eq11}
\end{equation}
Given $u\in X$, $t>0$, we define an operator on the space $Z:=C\left( \left[
0,t\right] ,L\left( X\right) \right) $, equipped with the supremum norm $%
\left\| \psi \right\| _Z=\sup\limits_{s\in \left[ 0,t\right] }\left\| \psi
\left( s\right) \right\| _{L\left( X\right) }$, by
\begin{equation}
\left( T\psi \right) \left( t\right) :=I+\int_0^tg^{\prime }\left( \tau
\left( s,u\right) \right) \psi \left( s\right) ds,\text{ }  \label{eq12}
\end{equation}
for all $\psi \in Z$, where $I$ stands for the identity map. Note that $%
\left( \ref{eq12}\right) $ is justified because the mapping $\varphi \left(
s\right) =I$ is constant for $s\in \left[ 0,t\right] $ and therefore $I\in
C\left( \left[ 0,t\right] ,L\left( X\right) \right) $. One can verify that $%
T $ maps $Z$ into itself, as the continuity of $T$ on $Z$ is ensured by the
boundedness of $s\rightarrow g^{\prime }\left( \tau \left( s,u\right)
\right) $ on $[0,t]$. Hence the equation $\left( \ref{eq11}\right) $ can be
recast as the following fixed-point problem
\begin{equation}
T\psi =\psi .  \label{eq13}
\end{equation}
For fixed $u\in X$, we define $A_{t,u}:=\bigcup\limits_{s\in \left[
0,t\right] }\left\{ \tau \left( s,u\right) \right\} $. Since $A_{t,u}$ is a
compact set and $g^{\prime }\left( \tau \left( s,u\right) \right) \in
L\left( X\right) $ for $s\in \left[ 0,t\right] $, there exists a constant $%
c_1\left( t\right) >0$ such that
\[
\sup\limits_{s\in \left[ 0,t\right] }\left\| g^{\prime }\left( \tau \left(
s,u\right) \right) \right\| _{L\left( X\right) }\leq c_1\left( t\right) .
\]

Notice that for any $\psi _1$,$\psi _2\in Z$,
\begin{eqnarray*}
\ \left\| T\psi _1-T\psi _2\right\| _Z &\leq &\int_0^t\sup\limits_{s\in
\left[ 0,t\right] }\left\| g^{\prime }\left( \tau \left( s,u\right) \right)
\left( \psi _1\left( s\right) -\psi _2\left( s\right) \right) \right\|
_{L\left( X\right) }ds \\
&\leq &c_1\left( t\right) t\cdot \left\| \psi _1-\psi _2\right\| _Z.
\end{eqnarray*}
For $t>0$ sufficiently small, the solvability of $\left( \ref{eq13}\right) $
is guaranteed by the contraction mapping principle; by the semigroup
property of the flow $\left( \ref{eq5}\right) $, it likewise holds for all $%
t>0$.

We now concentrate on the uniqueness of solutions of $\left( \ref{eq13}%
\right) $. Let $\widehat{\psi }_1$,$\widehat{\psi }_2\in Z$ be two solutions
of $\left( \ref{eq13}\right) $. In terms of the definition,
\begin{eqnarray*}
&&\ \left\| \widehat{\psi }_1\left( \xi \right) -\widehat{\psi }_2\left( \xi
\right) \right\| _{L\left( X\right) } \\
&\leq &\int_0^\xi \left\| g^{\prime }\left( \tau \left( s,u\right) \right)
\left( \widehat{\psi }_1\left( s\right) -\widehat{\psi }_2\left( s\right)
\right) \right\| _{L\left( X\right) }ds \\
&\leq &c_1\left( t\right) \int_0^\xi \left\| \widehat{\psi }_1\left(
s\right) -\widehat{\psi }_2\left( s\right) \right\| _{L\left( X\right) }ds
\end{eqnarray*}
and employing $\left( \ref{eq7}\right) $ we yield
\begin{equation}
\left\| \widehat{\psi }_1\left( \xi \right) -\widehat{\psi }_2\left( \xi
\right) \right\| _{L\left( X\right) }=0  \label{eq14}
\end{equation}
by taking $z\left( \xi \right) =\left\| \widehat{\psi }_1\left( \xi \right) -%
\widehat{\psi }_2\left( \xi \right) \right\| _{L\left( X\right) }$, $a=0$
and $b\left( s\right) =c_1\left( t\right) $ for $\xi \in \left[ 0,t\right] $%
. $\left( \ref{eq14}\right) $ shows that the solution of $\left( \ref{eq13}%
\right) $ is unique.

We now verify that $\tau \in C^2\left( \Bbb{R}^{+}\times X,X\right) $ and
actually just need to prove $\tau \left( t,\cdot \right) \in C^2\left(
X,X\right) $ for $t\in \Bbb{R}^{+}$ since it is unambiguous to check $\tau
\left( \cdot ,u\right) \in C^2\left( \Bbb{R}^{+},X\right) $ for $u\in X$.
Given $u\in X$, notice that for $\delta >0$, $v\in S:=\left\{ v\in X:\left\|
v\right\| _X=1\right\} $,
\begin{eqnarray}
&&\ \tau \left( \widetilde{t},u+\delta v\right) -\tau \left( \widetilde{t}%
,u\right)  \nonumber \\
\ &=&\delta v+\int_0^{\widetilde{t}}\left( g\left( \tau \left( s,u+\delta
v\right) \right) -g\left( \tau \left( s,u\right) \right) \right) ds
\nonumber \\
\ &=&\delta v+\int_0^{\widetilde{t}}\int_0^1g^{\prime }\left( \tau \left(
s,u\right) +\xi \left( \tau \left( s,u+\delta v\right) -\tau \left(
s,u\right) \right) \right)  \nonumber \\
&&\ \cdot \left( \tau \left( s,u+\delta v\right) -\tau \left( s,u\right)
\right) d\xi ds  \label{eq15}
\end{eqnarray}
for $\forall \widetilde{t}\in \left[ 0,t\right] $. Take $\delta _0>0$
sufficiently small and define
\[
\tau _{u,v}\left( \xi ,s,\delta \right) :=\left( 1-\xi \right) \tau \left(
s,u\right) +\xi \tau \left( s,u+\delta v\right)
\]
for fixed $\xi \in \left[ 0,1\right] $, $s\in \left[ 0,t\right] $, $\delta
\in \left[ 0,\delta _0\right] $.

Set $\varpi =\left( \xi ,s,\delta \right) \in Q:=\left[ 0,1\right] \times
\left[ 0,t\right] \times \left[ 0,\delta _0\right] $, $\tau _{u,v}\left(
\varpi \right) =\tau _{u,v}\left( \xi ,s,\delta \right) $. Due to
\[
\sup\limits_{\varpi \in Q}\left\| g^{\prime }\left( \tau _{u,v}\left( \varpi
\right) \right) \right\| _{L\left( X\right) }\leq c_2\left( t\right)
\]
for some $c_2\left( t\right) >0$, it follows from $\left( \ref{eq15}\right) $
that
\[
\left\| \tau \left( \widetilde{t},u+\delta v\right) -\tau \left( \widetilde{t%
},u\right) \right\| _X\leq \delta +c_2\left( t\right) \int_0^{\widetilde{t}%
}\left\| \tau \left( s,u+\delta v\right) -\tau \left( s,u\right) \right\|
_Xds.
\]
Take

\begin{eqnarray*}
z\left( \widetilde{t}\right) &=&\left\| \tau \left( \widetilde{t},u+\delta
v\right) -\tau \left( \widetilde{t},u\right) \right\| _X, \\
a &=&\delta , \\
b\left( s\right) &=&c_2\left( t\right) ,\text{\textrm{\ }}s\in \left[
0,t\right] .
\end{eqnarray*}
Again by $\left( \ref{eq7}\right) $, there exists a $c_3\left( t\right) >0$
such that
\begin{equation}
\left\| \tau \left( \widetilde{t},u+\delta v\right) -\tau \left( \widetilde{t%
},u\right) \right\| _X\leq c_3\left( t\right) \delta .  \label{eq16}
\end{equation}

Define $A_{u,v}\left( \varpi \right) :=g^{\prime }\left( \tau _{u,v}\left(
\varpi \right) \right) -g^{\prime }\left( \tau \left( s,u\right) \right) $.
Let $\varphi $ be a solution of $\left( \ref{eq13}\right) $. Hence
\begin{eqnarray*}
&&\tau \left( \widetilde{t},u+\delta v\right) -\tau \left( \widetilde{t}%
,u\right) -\varphi \left( \widetilde{t}\right) \delta v \\
&=&\int_0^{\widetilde{t}}\int_0^1g^{\prime }\left( \tau _{u,v}\left( \varpi
\right) \right) \left( \tau \left( s,u+\delta v\right) -\tau \left(
s,u\right) \right) d\xi ds \\
&&\ -\int_0^{\widetilde{t}}g^{\prime }\left( \tau \left( s,u\right) \right)
\varphi \left( s\right) \delta vds \\
&=&\int_0^{\widetilde{t}}g^{\prime }\left( \tau \left( s,u\right) \right)
\left( \tau \left( s,u+\delta v\right) -\tau \left( s,u\right) -\varphi
\left( s\right) \delta v\right) ds \\
&&+\int_0^{\widetilde{t}}\int_0^1A_{u,v}\left( \varpi \right) \left( \tau
\left( s,u+\delta v\right) -\tau \left( s,u\right) \right) d\xi ds,
\end{eqnarray*}
and this gets
\begin{eqnarray}
&&\ \ \left\| \tau \left( \widetilde{t},u+\delta v\right) -\tau \left(
\widetilde{t},u\right) -\varphi \left( \widetilde{t}\right) \delta v\right\|
_X  \nonumber \\
&\leq &\sup\limits_{\varpi \in Q}\left\| A_{u,v}\left( \varpi \right)
\right\| \cdot tc_3\left( t\right) \delta +c_1\left( t\right) \int_0^{%
\widetilde{t}}\left\| \ \tau \left( s,u+\delta v\right) -\tau \left(
s,u\right) -\varphi \left( s\right) \delta v\right\| _Xds.  \nonumber \\
&&\   \label{eq17}
\end{eqnarray}
Set
\[
\tau _{u,v}^{*}\left( \zeta \right) :=\tau _{u,v}^{*}\left( \mu ,\xi
,s,\delta \right) =\tau \left( s,u\right) +\mu \xi \left( \tau \left(
s,u+\delta v\right) -\tau \left( s,u\right) \right)
\]
for $\mu \in \left[ 0,1\right] $, $\zeta =\left( \mu ,\xi ,s,\delta \right)
\in Q^{*}=\left[ 0,1\right] \times \left[ 0,1\right] \times \left[
0,t\right] \times \left[ 0,\delta _0\right] $.

Observe that
\begin{eqnarray}
&&\ \ \ \left\| A_{u,v}\left( \varpi \right) \right\| _{L\left( X\right) }
\nonumber \\
\ &=&\left\| \int_0^1g^{\prime \prime }\left( \tau \left( s,u\right) +\mu
\xi \left( \tau \left( s,u+\delta v\right) -\tau \left( s,u\right) \right)
\right) \xi \left( \tau \left( s,u+\delta v\right) -\tau \left( s,u\right)
\right) d\mu \right\| _{L\left( X\right) }  \nonumber \\
\ &\leq &\sup\limits_{\zeta \in Q^{*}}\left\| g^{\prime \prime }\left( \tau
_{u,v}^{*}\left( \zeta \right) \right) \right\| _{L\left( X,L\left( X\right)
\right) }\cdot \left\| \tau \left( s,u+\delta v\right) -\tau \left(
s,u\right) \right\| _X  \label{eq18}
\end{eqnarray}
together with $g\in C^2\left( X,X\right) $, yields that there exists a $%
c_4\left( t\right) >0$, s.t.
\begin{equation}
\sup\limits_{\zeta \in Q^{*}}\left\| g^{\prime \prime }\left( \tau
_{u,v}^{*}\left( \zeta \right) \right) \right\| _{L\left( X,L\left( X\right)
\right) }\leq c_4\left( t\right) .  \label{eq19}
\end{equation}
Combining $\left( \ref{eq16}\right) $, $\left( \ref{eq17}\right) $, $\left(
\ref{eq18}\right) $ and $\left( \ref{eq19}\right) $ we obtain
\begin{eqnarray}
&&\ \ \ \ \left\| \tau \left( \widetilde{t},u+\delta v\right) -\tau \left(
\widetilde{t},u\right) -\varphi \left( \widetilde{t}\right) \delta v\right\|
_X  \nonumber \\
\ &\leq &c_5\left( t\right) \delta ^2+c_1\left( t\right) \int_0^{\widetilde{t%
}}\left\| \ \tau \left( s,u+\delta v\right) -\tau \left( s,u\right) -\varphi
\left( s\right) \delta v\right\| _Xds,  \label{eq20}
\end{eqnarray}
and by using $\left( \ref{eq7}\right) $ as before, $\left( \ref{eq20}\right)
$ gives rise to
\begin{equation}
\ \left\| \tau \left( \widetilde{t},u+\delta v\right) -\tau \left(
\widetilde{t},u\right) -\varphi \left( \widetilde{t}\right) \delta v\right\|
_X\leq c_5\left( t\right) \delta ^2\exp \left\{ c_1\left( t\right)
\widetilde{t}\right\}  \label{eq21}
\end{equation}
and subsequently produces
\begin{eqnarray}
\tau _u^{\prime }\left( \widetilde{t},u\right) &=&\varphi \left( \widetilde{t%
}\right) =I+\int_0^{\widetilde{t}}g^{\prime }\left( \tau \left( s,u\right)
\right) \varphi \left( s\right) ds  \nonumber \\
\ &=&I+\int_0^{\widetilde{t}}g^{\prime }\left( \tau \left( s,u\right)
\right) \tau _u^{\prime }\left( s,u\right) ds.  \label{eq22}
\end{eqnarray}
Given $\xi \in X\backslash \left\{ 0\right\} $, set $\phi _{u,\xi }\left(
\widetilde{t}\right) =\left\| \tau _u^{\prime }\left( \widetilde{t},u\right)
\xi \right\| _X$. As $\phi _{u,\xi }\left( 0\right) =\left\| \xi \right\| _X$%
, we have
\begin{equation}
\phi _{u,\xi }\left( \widetilde{t}\right) \leq \phi _{u,\xi }\left( 0\right)
+\max\limits_{s\in \left[ 0,\widetilde{t}\right] }\left\| g^{\prime }\left(
\tau \left( s,u\right) \right) \right\| _{L\left( X\right) }\cdot \int_0^{%
\widetilde{t}}\phi _{u,\xi }\left( s\right) ds  \label{eq23}
\end{equation}
for $\widetilde{t}>0$. By employing the first $``\leq "$ of $\left( \ref{eq9}%
\right) $,
\begin{equation}
\phi _{u,\xi }\left( \widetilde{t}\right) \geq \left\| \xi \right\| _X\exp
\left\{ -\max\limits_{s\in \left[ 0,\widetilde{t}\right] }\left\| g^{\prime
}\left( \tau \left( s,u\right) \right) \right\| _{L\left( X\right) }\cdot
\widetilde{t}\right\} >0.  \label{eq24}
\end{equation}
This concludes the proof of injectivity. We note that the argument relies on
the integral equation $\left( \ref{eq22}\right) $, which is satisfied by the
unique solution $\psi \left( t\right) =\tau _u^{\prime }\left( t,u\right) $%
(obtained as the fixed point of $T$), and from which the inequality $\left(
\ref{eq24}\right) $ is derived.

It remains to prove that $\tau \left( t,\cdot \right) \in C^2\left(
X,X\right) $. For fixed $t\in \Bbb{R}^{+}$, denote by a Banach space $%
\widetilde{Z}:=C\left( \left[ 0,t\right] ,L\left( X,L\left( X\right) \right)
\right) $, endowed with the supremum norm
\[
\Vert \Psi \Vert _{\widetilde{Z}}=\sup_{s\ \in [0,t]}\Vert \Psi (s)\Vert
_{L(X,L(X))}.
\]

We define an operator $\widetilde{T}:\widetilde{Z}\rightarrow \widetilde{Z}$
as follows:
\[
\left( \widetilde{T}\Psi \right) \left( t\right) v:=\int_0^t\left( g^{\prime
\prime }\left( \tau \left( s,u\right) \right) \tau _u^{\prime }\left(
s,u\right) v\right) \tau _u^{\prime }\left( s,u\right) ds+\int_0^tg^{\prime
}\left( \tau \left( s,u\right) \right) \Psi \left( s\right) vds,
\]
for all $v\in X$ and $\Psi \in \widetilde{Z}$.

Plainly $\widetilde{T}$ is well-defined on $\widetilde{Z}$ by virtue of the
facts that $g^{\prime \prime }\left( \tau \left( s,u\right) \right) \tau
_u^{\prime }\left( s,u\right) \in L\left( X,L\left( X\right) \right) $, and $%
g^{\prime }\left( \tau \left( s,u\right) \right) \in L\left( X\right) $ for $%
s\in \left[ 0,t\right] $. Consider the following fixed point problem
\begin{equation}
\widetilde{T}\Psi =\Psi .  \label{eq25}
\end{equation}
Given $\Psi _1$,$\Psi _2\in \widetilde{Z}$, since $A_{t,u}$ is a compact set
and $g^{\prime }\left( \tau \left( s,u\right) \right) \in L\left( X\right) $
for $s\in \left[ 0,t\right] $, there exists a $c_6\left( t\right) >0$ such
that
\begin{eqnarray}
\ \left\| \widetilde{T}\Psi _1-\widetilde{T}\Psi _2\right\| _{\widetilde{Z}}
&=&\sup\limits_{v\in X,\text{ }\left\| v\right\| =1}\left\| \left(
\widetilde{T}\Psi _1\right) v-\left( \widetilde{T}\Psi _2\right) v\right\|
_{L\left( X\right) }  \nonumber \\
&\leq &\int_0^t\sup\limits_{s\in \left[ 0,t\right] }\left\| g^{\prime
}\left( \tau \left( s,u\right) \right) \left( \Psi _1\left( s\right) -\Psi
_2\left( s\right) \right) \right\| _{L\left( X,L\left( X\right) \right) }ds
\nonumber \\
&\leq &c_6\left( t\right) t\left\| \Psi _1-\Psi _2\right\| _{\widetilde{Z}}.
\label{eq26}
\end{eqnarray}
and also for the solutions of $\left( \ref{eq25}\right) $, denoted by $%
\widehat{\Psi }_1$,$\widehat{\Psi }_2\in \widetilde{Z}$,
\begin{equation}
\left\| \widehat{\Psi }_1\left( \xi \right) -\widehat{\Psi }_2\left( \xi
\right) \right\| _{L\left( X,L\left( X\right) \right) }\leq c_6\left(
t\right) \int_0^\xi \left\| \widehat{\Psi }_1\left( s\right) -\widehat{\Psi }%
_2\left( s\right) \right\| _{L\left( X,L\left( X\right) \right) }ds
\label{eq27}
\end{equation}
for $\xi \in \left[ 0,t\right] $, where $\left\| \Psi \right\| _{\widetilde{Z%
}}=\sup\limits_{s\in \left[ 0,t\right] }\left\| \Psi \left( s\right)
\right\| _{L\left( X,L\left( X\right) \right) }.$

Analogously to the preceding argument, the solvability and uniqueness of $%
\left( \ref{eq25}\right) $ can be readily established. Invoking $\left( \ref
{eq7}\right) $, one deduces that there exists a unique solution $\widehat{%
\Psi }\left( s\right) \in \widetilde{Z}$ solving $\left( \ref{eq25}\right) $.

We infer from $\left( \ref{eq22}\right) $ that
\begin{eqnarray}
&&\ \ \tau _u^{\prime }\left( t,u+\delta v\right) -\tau _u^{\prime }\left(
t,u\right)  \nonumber \\
\ &=&\int_0^tg^{\prime }\left( \tau \left( s,u+\delta v\right) \right) \tau
_u^{\prime }\left( s,u+\delta v\right) -g^{\prime }\left( \tau \left(
s,u\right) \right) \tau _u^{\prime }\left( s,u\right) ds  \nonumber \\
\ &=&\int_0^t\int_0^1g^{\prime \prime }\left( \tau \left( s,u\right) +\xi
\left( \tau \left( s,u+\delta v\right) -\tau \left( s,u\right) \right)
\right) \left( \tau \left( s,u+\delta v\right) -\tau \left( s,u\right)
\right)  \nonumber \\
&&\ \ \cdot \tau _u^{\prime }\left( s,u+\delta v\right) d\xi
ds+\int_0^tg^{\prime }\left( \tau \left( s,u\right) \right) \left[ \tau
_u^{\prime }\left( s,u+\delta v\right) -\tau _u^{\prime }\left( s,u\right)
\right] ds.  \label{eq28}
\end{eqnarray}
The assumption $g\in C^2\left( \Bbb{R}^{+}\times X,X\right) $, together with
$\left( \ref{eq7}\right) $, $\left( \ref{eq16}\right) $ and $\left( \ref
{eq28}\right) $, yields that there exists a $c_7\left( t\right) >0$, s.t.
\begin{equation}
\left\| \tau _u^{\prime }\left( t,u+\delta v\right) -\tau _u^{\prime }\left(
t,u\right) \right\| _{L\left( X\right) }\leq c_7\left( t\right) \delta .
\label{eq29}
\end{equation}

Since $\widehat{\Psi }\left( t\right) $ solves $\left( \ref{eq25}\right) $,
applying $\left( \ref{eq28}\right) $ we derive
\begin{eqnarray}
&&\ \tau _u^{\prime }\left( t,u+\delta v\right) -\tau _u^{\prime }\left(
t,u\right) -\widehat{\Psi }\left( t\right) \delta v  \nonumber \\
\ &=&\int_0^t\int_0^1g^{\prime \prime }\left( \tau \left( s,u\right) +\xi
\left( \tau \left( s,u+\delta v\right) -\tau \left( s,u\right) \right)
\right) \left( \tau \left( s,u+\delta v\right) -\tau \left( s,u\right)
\right)  \nonumber \\
&&\ \cdot \tau _u^{\prime }\left( s,u+\delta v\right) d\xi
ds+\int_0^tg^{\prime }\left( \tau \left( s,u\right) \right) \left[ \tau
_u^{\prime }\left( s,u+\delta v\right) -\tau _u^{\prime }\left( s,u\right)
\right] ds  \nonumber \\
&&\ -\int_0^tg^{\prime \prime }\left( \tau \left( s,u\right) \right) \tau
_u^{\prime }\left( s,u\right) \delta v\tau _u^{\prime }\left( s,u\right)
ds-\int_0^tg^{\prime }\left( \tau \left( s,u\right) \right) \widehat{\Psi }%
\left( s\right) \delta vds  \nonumber \\
\ &=&\int_0^tg^{\prime }\left( \tau \left( s,u\right) \right) \left( \tau
_u^{\prime }\left( s,u+\delta v\right) -\tau _u^{\prime }\left( s,u\right) -%
\widehat{\Psi }\left( s\right) \delta v\right) ds  \nonumber \\
&&\ +\int_0^t\int_0^1g^{\prime \prime }\left( \tau \left( s,u\right) +\xi
\left( \tau \left( s,u+\delta v\right) -\tau \left( s,u\right) \right)
\right) \left( \tau \left( s,u+\delta v\right) -\tau \left( s,u\right)
\right)  \nonumber \\
&&\ \cdot \tau _u^{\prime }\left( s,u+\delta v\right) d\xi
ds-\int_0^tg^{\prime \prime }\left( \tau \left( s,u\right) \right) \tau
_u^{\prime }\left( s,u\right) \delta v\tau _u^{\prime }\left( s,u\right) ds.
\label{eq30}
\end{eqnarray}
By $\left( \ref{eq7}\right) $,
\begin{equation}
\left\| \tau _u^{\prime }\left( t,u+\delta v\right) -\tau _u^{\prime }\left(
t,u\right) -\widehat{\Psi }\left( t\right) \delta v\right\| _{L\left(
X\right) }\leq ta_t\left( \delta v\right) \cdot \exp \int_0^tB_1\left(
s\right) ds,  \label{eq31}
\end{equation}
where
\begin{eqnarray*}
B_1\left( t\right) &=&\sup\limits_{0\leq s\leq t}\left\| g^{\prime }\left(
\tau \left( s,u\right) \right) \right\| _{L\left( X\right) }, \\
a_t\left( \delta v\right) &=&\sup\limits_{0\leq s\leq t}\left\|
\int_0^1C\left( s,\xi ,u,\delta v\right) -g^{\prime \prime }\left( \tau
\left( s,u\right) \right) \tau _u^{\prime }\left( s,u\right) \delta v\tau
_u^{\prime }\left( s,u\right) d\xi \right\| _{L\left( X\right) },
\end{eqnarray*}
and
\begin{eqnarray*}
\ C\left( s,\xi ,u,\delta v\right) &=&g^{\prime \prime }\left( \tau \left(
s,u\right) +\xi \left( \tau \left( s,u+\delta v\right) -\tau \left(
s,u\right) \right) \right) \\
&&\cdot \left( \tau \left( s,u+\delta v\right) -\tau \left( s,u\right)
\right) \tau _u^{\prime }\left( s,u+\delta v\right) .
\end{eqnarray*}

The assumption that $g\in C^2\left( \Bbb{R}^{+}\times X,X\right) $, combined
with $\left( \ref{eq15}\right) $, allows us to apply Lebesgue dominated
convergence theorem to obtain that
\begin{equation}
a_t\left( \delta v\right) =o\left( \delta \right)  \label{eq32}
\end{equation}
as $s\in \left[ 0,t\right] $ and $\delta >0$ sufficiently small. By $\left(
\ref{eq31}\right) $ and $\left( \ref{eq32}\right) $ we get
\[
\tau _u^{\prime \prime }\left( t,u\right) v=\widehat{\Psi }\left( t\right)
v=\int_0^tg^{\prime \prime }\left( \tau \left( s,u\right) \right) \tau
_u^{\prime }\left( s,u\right) v\tau _u^{\prime }\left( s,u\right) +g^{\prime
}\left( \tau \left( s,u\right) \right) \widehat{\Psi }\left( s\right) vds.
\]
The $C^2$ smoothness of solutions is confirmed. That's the statement of the
first assertion of the lemma. The last assertion follows directly from
noting that
\begin{eqnarray*}
\left\| \frac{d\tau \left( t,u\right) }{dt}\right\| _X &\leq &\left\|
g\left( u\right) \right\| _X+\left\| \int_0^t\frac{dg\left( \tau \left(
s,u\right) \right) }{ds}ds\right\| _X \\
&=&\left\| g\left( u\right) \right\| _X+\left\| \int_0^tg^{\prime }\left(
\tau \left( s,u\right) \right) \frac{d\tau \left( s,u\right) }{ds}ds\right\|
_X \\
&\leq &\left\| g\left( u\right) \right\| _X+C_t\int_0^t\left\| \frac{d\tau
\left( s,u\right) }{ds}\right\| _Xds
\end{eqnarray*}
yields
\begin{equation}
\left\| g\left( u\right) \right\| _X\exp \left\{ -C_t\cdot t\right\} \leq
\left\| \frac{d\tau \left( t,u\right) }{dt}\right\| _X\leq \left\| g\left(
u\right) \right\| _X\exp \left\{ C_t\cdot t\right\}  \label{eq33}
\end{equation}
by $\left( \ref{eq9}\right) $, where $C_t:=$ $\max\limits_{s\in \left[
0,t\right] }\left\| g^{\prime }\left( \tau \left( s,u\right) \right)
\right\| _X$. $\left( \ref{eq33}\right) $ also shows that $g\left( \tau
\left( t,u\right) \right) =0$ for $\forall t>0$ if $g\left( u\right) =0$. We
therefore conclude the proof.\qed\vskip 5pt

\section{The smoothing of $f\left( x,t\right) $}

The requirement of $C^2$ smoothness on the nonlinear term $f\left(
x,t\right) $ is too restrictive for the original problem. Typically, we only
assume $f$ is $C^1$(see condition $\left( f_1^{\prime }\right) $). For this
reason, a direct proof of Theorem \ref{thm1.1} is infeasible, and we must
first mollify the nonlinearity. We begin by recalling the standard mollifier
in $\Bbb{R}^M$(see \cite{Br}):

\[
\rho \left( t\right) =\left\{
\begin{array}{cc}
e^{\frac 1{\left| t\right| ^2-1}}, & \left| t\right| <1, \\
0, & \left| t\right| \geq 1,
\end{array}
\right.
\]
then $\rho \left( t\right) \in C_c^\infty \left( \Bbb{R}^M\right) $. Let $%
\rho _m\left( t\right) =cm^M\rho \left( mt\right) $, $m\in \Bbb{N}$, $%
c=\left( \int_{\Bbb{R}^M}\rho \right) ^{-1}$. Define
\[
g_m\left( t\right) :=\rho _m*g=\int_{\Bbb{R}^M}\rho _m\left( t-s\right)
g\left( s\right) ds.
\]

\begin{proposition}
\label{prop3.1}(see \cite{Br}) $\left( i\right) $ If $g\in L_{loc}^1\left(
\Bbb{R}^M\right) $, then $\rho _m*g$ $\in C^\infty \left( \Bbb{R}^M\right) $%
, and $D^\alpha \left( \rho _m*g\right) =\left( D^\alpha \rho _m\right) *g$,
where $D^\alpha \rho _m=\frac{\partial ^{\alpha _1}}{\partial t_1^{\alpha _1}%
}\cdots \frac{\partial ^{\alpha _M}}{\partial t_M^{\alpha _M}}\rho _m$, $%
\left| \alpha \right| =\alpha _1+\cdots +\alpha _M$;

$\left( ii\right) $ If $g\in C\left( \Bbb{R}^M\right) $, then $\rho _m*g$
uniformly converges to $g$ in any compact set of $\Bbb{R}^M$;

$\left( iii\right) $ If $g\in L^p\left( \Bbb{R}^M\right) $, $1\leq p<\infty $%
, then $\rho _m*g$ converges to $g$ in $L^p\left( \Bbb{R}^M\right) $.
\end{proposition}

Let $M=1$ and denote

\[
\text{\textrm{\ }}f_m\left( x,t\right) :=\left( \rho _m*f\right) \left(
x,t\right) =\int_{-\infty }^\infty \rho _m\left( t-s\right) f\left(
x,s\right) ds.
\]
In view of the definition,
\[
\text{\textrm{\ }}f_m\left( x,t\right) =\int_{-\frac 1m}^{\frac 1m}\rho
_m\left( \tau \right) f\left( x,t-\tau \right) d\tau
\]
and Proposition \ref{prop3.1} accordingly yields
\begin{equation}
f_m\left( x,t\right) \rightarrow \text{\textrm{\ }}f\left( x,t\right)
\label{eq34}
\end{equation}
uniformly in any compact set of $\overline{\Omega }\times \Bbb{R}$. In
addition,
\begin{equation}
\left( f_m\right) _t^{\left( i\right) }\left( x,t\right) =\int_{-\infty
}^\infty \rho _m^{\left( i\right) }\left( \tau \right) f\left( x,t-\tau
\right) d\tau .  \label{eq35}
\end{equation}

\begin{lemma}
\label{lem3.2} For $m\in \Bbb{N}$ large, if $f\left( x,t\right) $ satisfies
the hypothesis $\left( f_2\right) $, $f_m\left( x,t\right) $ does as well.
\end{lemma}

\textit{Proof}. Take $\widetilde{\theta }\in \left( 2,\theta \right) $. It
suffices to show that $\exists M_1>0$, $\exists N_0^{\left( 1\right) }\in
\Bbb{N}$, s.t. for $x\in \overline{\Omega }$, and $m\in \Bbb{N}$, $m\geq
N_0^{\left( 1\right) }$,
\begin{equation}
0<\widetilde{\theta }F_m\left( x,t\right) \leq tf_m\left( x,t\right)
\label{eq36}
\end{equation}
as $\left| t\right| >M_1$, where $F_m\left( x,t\right) =\int_0^tf_m\left(
x,s\right) ds$. Indeed, a direct computation gives rise to
\begin{eqnarray*}
&&\ \ \ \ \ \ \ \ \ tf_m\left( x,t\right) -\widetilde{\theta }%
\int_0^tf_m\left( x,s\right) ds \\
\ \ &=&\int_{-\frac 1m}^{\frac 1m}\rho _m\left( \tau \right) t\cdot f\left(
x,t-\tau \right) d\tau -\widetilde{\theta }\int_0^t\int_{-\frac 1m}^{\frac
1m}\rho _m\left( \tau \right) \cdot f\left( x,s-\tau \right) d\tau ds \\
\ &=&\int_{-\frac 1m}^{\frac 1m}\rho _m\left( \tau \right) \left( t\cdot
f\left( x,t-\tau \right) -\widetilde{\theta }\int_0^tf\left( x,s-\tau
\right) ds\right) d\tau ,
\end{eqnarray*}
so we merely need to verify that there exist constants $\beta >0$, $M_1>0$,
such that
\[
t\cdot f\left( x,t-\tau \right) -\widetilde{\theta }\int_0^tf\left( x,s-\tau
\right) ds>\beta
\]
as $\left| t\right| >M_1$. Set $\xi =t-\tau $. In terms of the hypothesis $%
\left( f_2\right) $, there exist constants $\alpha _1>0$, $\alpha _2>0$, $%
\beta >0$, such that
\begin{eqnarray*}
&&\ \ t\cdot f\left( x,\xi \right) -\widetilde{\theta }\int_0^tf\left(
x,s-\tau \right) ds \\
\ &\geq &\left( \frac{\xi +\tau }\xi \theta -\widetilde{\theta }\right)
F\left( x,\xi \right) +\widetilde{\theta }F\left( x,-\tau \right) \\
\ &\geq &\alpha _1\left| \xi \right| ^\theta -\alpha _2>\beta
\end{eqnarray*}
as $\left| t\right| \geq M_1:=M+1$, $m\in \Bbb{N}$ large. We thereby arrive
at the conclusion.\qed\vskip 5pt

Consider the equation
\begin{equation}
\left\{
\begin{array}{cc}
-\Delta u=f_m\left( x,u\right) , & x\in \Omega , \\
u=0, & x\in \partial \Omega ,
\end{array}
\right.  \label{eq37}
\end{equation}
and define the corresponding energy functionals
\[
J_m\left( u\right) =\frac 12\int_\Omega \left| \nabla u\right|
^2dx-\int_\Omega F_m\left( x,u\right) dx,
\]
where $F_m\left( x,t\right) =\int_0^tf_m\left( x,s\right) ds.$

Applying Lemma \ref{lem3.2} it is straightforward to check that $\Phi
_m\left( u\right) :=\int_\Omega F_m\left( x,u\right) dx$ is a superquadratic
function defined on $Y:=L^{\widetilde{\theta }}\left( \Omega \right) $, $%
\widetilde{\theta }$ $\in \left( 2,\theta \right) $. Let $\left\{ \varphi
_k\right\} _{k=1}^\infty $ be the orthonormal eigenvectors of the linear
eigenvalue problem $\left( \ref{eq93}\right) $. The spectral theorem
guarantees that $\left\{ \varphi _k\right\} _{k=1}^\infty $ forms a complete
orthonormal basis of $X$. For $n\in \Bbb{N}$, define the finite-dimensional
subspace $X_n:=$span$\left\{ \varphi _1,\cdots ,\varphi _n\right\} $. The
next lemma concerns the smoothness of the restricted functional:

\begin{lemma}
\label{lem3.3} Under assumption $\left( f_1^{\prime }\right) $, the
restriction $\Phi _{m,n}:=\Phi _m|_{X_n}$ is infinitely differentiable,
i.e., $\Phi _{m,n}\in C^\infty \left( X_n,R\right) $, where $\Phi
_{m,n}\left( u\right) :=\int_\Omega F_m\left( x,u\right) dx$ for $u\in X_n$.
\end{lemma}

\textit{Proof}. Observe that for $x\in \overline{\Omega }$, $u_0\in X_n$, $%
h\in S^{n-1}:=\left\{ x\in X_n:\left\| x\right\| _X=1\right\} $,
\begin{eqnarray*}
&&\ \ \ F_m\left( x,u_0\left( x\right) +sh\left( x\right) \right) -F_m\left(
x,u_0\left( x\right) \right) \\
\ &=&sf_m\left( x,u_0\left( x\right) +\xi \left( s,x\right) sh\left(
x\right) \right) h\left( x\right) ,\text{\textrm{\ } }0<\xi \left(
s,x\right) <1,
\end{eqnarray*}
via mean value theorem. Hence, $\forall \varepsilon >0$, $\exists \delta >0$%
, $0<\left| s\right| <\delta $,
\begin{eqnarray*}
&&\ \ \ \left| \frac 1s\left[ F_m\left( x,u_0\left( x\right) +sh\left(
x\right) \right) -F_m\left( x,u_0\left( x\right) \right) \right] -f_m\left(
x,u_0\left( x\right) \right) h\left( x\right) \right| \\
\ &=&\left| \left[ f_m\left( x,u_0\left( x\right) +\xi \left( s,x\right)
sh\left( x\right) \right) -f_m\left( x,u_0\left( x\right) \right) \right]
h\left( x\right) \right| \\
\ &\leq &\varepsilon \cdot \left| h\left( x\right) \right|
\end{eqnarray*}
as $h\left( x\right) \in C^\infty \left( \overline{\Omega }\right) $.
Consequently
\[
\frac 1s\left( F_m\left( x,u_0\left( x\right) +sh\left( x\right) \right)
-F_m\left( x,u_0\left( x\right) \right) \right) \rightarrow f_m\left(
x,u_0\left( x\right) \right) h\left( x\right)
\]
uniformly on $x\in \overline{\Omega }$ as $s\rightarrow 0$. Using Lebesgue
dominated convergence theorem,

\begin{eqnarray*}
&&\lim_{s\rightarrow 0}\frac 1s\int_\Omega F_m\left( x,u_0\left( x\right)
+sh\left( x\right) \right) -F_m\left( x,u_0\left( x\right) \right) dx \\
&=&\int_\Omega f_m\left( x,u_0\left( x\right) \right) h\left( x\right) dx,
\end{eqnarray*}
and this yields
\[
\Phi _{m,n}^{\prime }\left( u_0\right) =P_n\left( -\Delta \right)
^{-1}f_m\left( x,u_0\right) ,
\]
where $P_n$ denotes the orthogonal projection from $X$ to $X_n$.

Aiming at proving Lemma \ref{lem3.3}, we argue by induction. Without loss of
generality, for $u_0\in X_n$, $h_i\in S^{n-1}$, $i=1$,$2$,$\cdots $,$k$, we
have

\[
\Phi _{m,n}^{\left( k\right) }\left( x,u_0\right) \left( h_1,\cdots
,h_k\right) =\int_\Omega \left( f_m\right) _t^{\left( k-1\right) }\left(
x,u_0\right) h_1\cdots h_kdx,
\]
where
\[
\left( f_m\right) _t^{\left( i\right) }\left( x,t\right) =D^{\left( i\right)
}\rho _m*f\left( x,t\right) =\int_{-\frac 1m}^{\frac 1m}\rho _m^{\left(
i\right) }\left( \tau \right) f\left( x,t-\tau \right) d\tau .
\]

Observe that for $m\in \Bbb{N}$ large, $\tau \in \left[ -\frac 1m,\frac
1m\right] $, $\exists \widetilde{c}>0$, $d_{i,m}>0$ such that
\[
\left| f\left( x,t-\tau \right) \right| <\widetilde{c}\left( 1+\left|
t\right| ^{p-1}\right)
\]
and
\begin{eqnarray*}
\left| \left( f_m\right) _t^{\left( i\right) }\left( x,t\right) \right|
&\leq &\widetilde{c}\left( 1+\left| t\right| ^{p-1}\right) \int_{-\frac
1m}^{\frac 1m}\left| \rho _m^{\left( i\right) }\left( \tau \right) \right|
d\tau \\
&\leq &d_{i,m}\left( 1+\left| t\right| ^{p-1}\right)
\end{eqnarray*}
via $\left( \ref{eq35}\right) $, where $\widetilde{c}$ is independent of $%
\tau $ and $\rho _m^{\left( i\right) }\in C^\infty \left[ -\frac 1m,\frac
1m\right] $. Due to $u_0$,$h_{k+1}\in C_0^\infty \left( \overline{\Omega }%
\right) $, an analogous argument obtains
\begin{eqnarray*}
&&\ \lim_{s\rightarrow 0}\frac 1s\int_\Omega \left( \left( f_m\right)
_t^{\left( k-1\right) }\left( x,u_0+sh_{k+1}\right) h_1\cdots h_k-\left(
f_m\right) _t^{\left( k-1\right) }\left( x,u_0\right) h_1\cdots h_k\right) dx
\\
\ &=&\int_\Omega \left( f_m\right) _t^{\left( k\right) }\left( x,u_0\right)
h_1\cdots h_{k+1}dx.
\end{eqnarray*}
The assertion follows.\qed\vskip 5pt

\section{Proofs of Theorems 1.1 and 1.2}

Let $X$ be a Hilbert space endowed with the norm $\left\| \cdot \right\| _X$%
. We first take a look back at the concept of Palais-Smale condition. A
function $I$ is said to satisfy the Palais-Smale((PS) in short) condition on
a Hilbert space $X$, if any sequence $\left\{ x_i\right\} \subset X$ with $%
\left| I\left( x_i\right) \right| $ bounded and $\left\| I^{\prime }\left(
x_i\right) \right\| _X\rightarrow 0$ contains a convergent subsequence $%
\left\{ x_{i_k}\right\} $ in $X$. Let $\left\{ \varphi _i\right\}
_{i=1}^\infty $ be a complete orthonormal basis of $X$, $X_n=$span$\left\{
\varphi _1,\cdots ,\varphi _n\right\} $ and $I_n=I|_{X_n}$. Consider a
sequence $\left\{ X_n\right\} _{n=1}^\infty $ such that $X=\overline{%
\bigcup\limits_{n=1}^\infty X_n}$. We say that $I$ satisfies the (PS)$^{%
\text{*}}$ condition with respect to $\left\{ X_n\right\} _{n=1}^\infty $ if
every sequence $\left\{ u_n\right\} _{n=1}^\infty $, $u_n\subset X_n$, such
that
\[
\sup\limits_{n\in \Bbb{N}}I_n\left( u_n\right) <\infty ,\text{\textrm{\ }}%
I_n^{\prime }\left( u_n\right) \rightarrow 0\text{\textrm{\ as }}%
n\rightarrow \infty ,
\]
contains a subsequence convergent to a critical point of $I$(see \cite{Ba},
\cite{LW}).

Henceforth, unless otherwise specified, we will take $X=H_0^1\left( \Omega
\right) $. Aiming at verifying Theorem \ref{thm1.1}, the following
consequence stemming from \cite{Ghou} is adapted to our needs and we restate
the proof for the convenience of later discussion.

\begin{proposition}
\label{prop4.1} Let $E$ be a Hilbert space. If $K$ is a compact subset of $E$%
, then for every $\mu >0$ there exists a $C^\infty $ function $%
l:E\rightarrow \left[ 0,1\right] $ with all its derivatives bounded such
that
\[
l\left( x\right) =1\text{\textrm{,}}\mathrm{\ }\text{\textrm{for all }}x\in
K^\mu \text{\textrm{; }}l\left( x\right) =0\text{\textrm{,\ for all }}x\in
E\backslash K^{2\mu }\mathrm{,}
\]
where $K^\mu =\left\{ x\in E:\inf\limits_{z\in K}\left\| x-z\right\| _E<\mu
\right\} $.
\end{proposition}

\textit{Proof. }Fix a finite number of points $x_1$, $x_2$, $\cdots $, $%
x_m\in K$ such that $K\subset \bigcup\limits_{i=1}^mB\left( x_i,\frac \mu
2\right) $. We get $C^\infty $-functions $l_i$ on $E$ with bounded
derivatives such that
\[
l_i\left( x\right) =1\text{\textrm{\ as } }x\in B\left( x_i,\frac 32\mu
\right) \text{\textrm{; }}l_i\left( x\right) =0\text{\textrm{\ as } }x\notin
B\left( x_i,2\mu \right) \text{\textrm{,}}
\]
by simply taking a real function $\alpha :\Bbb{R\rightarrow }\left[
0,1\right] $ such that
\[
\alpha \left( t\right) =1\text{\textrm{\ } \textrm{if } }t\leq \frac 94\text{%
\textrm{; }}\alpha \left( t\right) =0\text{\textrm{\ if } }t\geq 4\text{%
\textrm{,}}
\]
and setting $l_i\left( x\right) =\alpha \left( \frac{\left\| x-x_i\right\|
_E^2}{\mu ^2}\right) $. Then we take $\beta :\Bbb{R}\rightarrow \left[
0,1\right] $ such that
\[
\beta \left( t\right) =0\mathrm{\ }\text{\textrm{if } }t\leq 0\text{\textrm{%
; }}\beta \left( t\right) =1\text{\textrm{\ if } }t\geq 1\text{\textrm{.}}
\]
The function $l:E\rightarrow \Bbb{R}$ defined by $l\left( u\right) =\beta
\left( \sum\limits_{i=1}^ml_i\left( u\right) \right) $ is the required
function.

Before entering the proof of Theorem \ref{thm1.1}, we recall a well-known
result (see \cite{Willem}):

\begin{proposition}
\label{prop4.2} Assume that $\left| \Omega \right| <\infty $, $1\leq p$, $%
r<\infty $, $g\in C\left( \overline{\Omega }\times \Bbb{R},\Bbb{R}\right) $
and
\begin{equation}
\left| g\left( x,u\right) \right| \leq C\left( 1+\left| u\right|
^{p/r}\right) \text{.}  \label{eq38}
\end{equation}
Then, for every $u\in L^p\left( \Omega \right) $, $g\left( \cdot ,u\right)
\in L^r\left( \Omega \right) $ and the operator
\[
A:L^p\left( \Omega \right) \rightarrow L^r\left( \Omega \right)
:u\longmapsto g\left( x,u\right)
\]
is continuous.
\end{proposition}

Divide the proof of Theorem \ref{thm1.1} $\left( i\right) $ into six steps:

Step 1\textbf{. }Observe that $\left( f_3\right) $ yields that for $%
\varepsilon \in \left( 0,\lambda _1\right) $, $\exists C_\varepsilon >0$,
\begin{equation}
\left| f\left( x,t\right) \right| \leq \varepsilon \left| t\right|
+C_\varepsilon \left| t\right| ^{p-1},\text{\textrm{\ }}\forall t\in \Bbb{R}.
\label{eq39}
\end{equation}
Hence for $u\in X$, $\left\| u\right\| _X=r>0$ sufficiently small, $\exists
\delta ^{*}>0$ s.t.
\begin{eqnarray}
J\left( u\right) &=&\frac 12\left\| u\right\| _X^2-\int_\Omega F\left(
x,u\right) dx  \nonumber \\
\ &\geq &\frac 12\left( 1-\frac \varepsilon {\lambda _1}\right) \left\|
u\right\| _X^2-\frac{C_\varepsilon }p\int_\Omega \left| u\right| ^pdx
\nonumber \\
\ &\geq &\left\| u\right\| _X^2\cdot \left[ \frac 12\left( 1-\frac
\varepsilon {\lambda _1}\right) -\frac{C_\varepsilon C_1}p\left\| u\right\|
_X^{p-2}\right]  \nonumber \\
\ &>&\delta ^{*}>0  \label{eq40}
\end{eqnarray}
with $C_1>0$. And also notice that $\left( f_2\right) $ gets
\begin{equation}
\int_\Omega f\left( x,u\right) u\geq C_2\int_\Omega \left| u\right| ^\theta
-C_3,  \label{eq41}
\end{equation}
for some $C_2$,$C_3>0$. Set $u=R\widehat{u}$, $\widehat{u}=\frac u{\left\|
u\right\| _X}$. Let $X_n=$span$\left\{ \varphi _1,\cdots ,\varphi _n\right\}
$, where $\varphi _i$ is the eigenfunction associated with the $i$-th
eigenvalue of the linear eigenvalue problem $\left( \ref{eq93}\right) $, $%
i=1 $,$2\cdots $,$n$. Therefore for $u\in X_n$, $\left\| u\right\| _X=R>0$,
\begin{eqnarray}
J\left( u\right) &=&\frac{R^2}2-\int_\Omega F\left( x,u\right) dx  \nonumber
\\
\ &\leq &\frac{R^2}2-C_4R^\theta \int_\Omega \left| \widehat{u}\right|
^\theta +C_5  \label{eq42}
\end{eqnarray}
as $R>0$ large. Given $k\in \Bbb{N}$, $k<n$. Using $\left( \ref{eq40}\right)
$ we can take $\rho _0>0$, s.t.
\[
d_0:=\inf\limits_{u\in X,\text{ }\left\| u\right\| _X=\rho _0}J\left(
u\right) >0.
\]
Take $r=r_k>0$ suitably small and $R=R_k>0$ large enough, s.t., $\rho _0\in
\left( r_k,R_k\right) $, and
\begin{equation}
\sup\limits_{u\in \overline{B_k\left( 0,r_k\right) }\cup \partial B_k\left(
0,R_k\right) }J\left( u\right) <\frac{d_0}2,  \label{eq43}
\end{equation}
where
\begin{eqnarray*}
B_k\left( 0,r_k\right) &:&=\left\{ u\in X_k:\left\| u\right\| _X<r_k\right\}
; \\
\partial B_k\left( 0,R_k\right) &:&=\left\{ u\in X_k:\left\| u\right\|
_X=R_k\right\} .
\end{eqnarray*}

Set $K=\overline{B_k\left( 0,r_k\right) }\cup \partial B_k\left(
0,R_k\right) $. Pick $\mu >0$ small enough such that
\[
\sup\limits_{x\in K_{X_n}^{2\mu }}J_n\left( x\right) <\frac{d_0}2,
\]
where $K_{X_n}^{2\mu }:=\left\{ x\in X_n:\inf\limits_{z\in K}\left\|
x-z\right\| _X<2\mu \right\} $ and $J_n:=J|_{X_n}$. Also denote by $J_{m,n}$
the restriction of the energy functional $J_m$ to $X_n$; by Lemma \ref
{lem3.3}, it is of class $C^\infty \left( X_n,R\right) $.

Invoking Proposition \ref{prop4.1}, there exists a $C^\infty $ function $%
l:X_n\rightarrow \left[ 0,1\right] $ satisfying
\[
l\left( x\right) =1\text{ for }x\in K_{X_n}^\mu ,l\left( x\right) =0\text{
for }x\in X_n\backslash K_{X_n}^{2\mu }.
\]
Now consider the following flow:
\begin{equation}
\left\{
\begin{array}{c}
\frac{d\tau \left( t,u\right) }{dt}=-\left( 1-l\left( \tau \left( t,u\right)
\right) \right) J_{m,n}^{\prime }\left( \tau \left( t,u\right) \right) ,%
\text{ }\left( t,u\right) \in \Bbb{R}^{+}\times X_n, \\
\tau \left( 0,u\right) =u\in X_n,
\end{array}
\right.  \label{eq44}
\end{equation}
where $J_{m,n}^{\prime }\left( u\right) :=P_nJ_m^{\prime }\left( u\right) $
and $P_n$ stands for the orthogonal projection from $X$ to $X_n$. Lemma \ref
{lem2.1} shows that $\tau \in C^2\left( \left[ 0,\infty \right) \times
X_n,X_n\right) $. Set $K_{J_{m,n}}:=\left\{ u\in X_n:J_{m,n}^{\prime }\left(
u\right) =0\right\} $. For each initial point $u\in X_n$, denote by $T_u$
the arrivial time at $K_{J_{m,n}}$ under the flow $\left( \ref{eq44}\right) $%
. By Lemma \ref{lem2.1}, $T_u=\infty $ if $u\notin K_{J_{m,n}}$, and $\tau
\left( t,u\right) =u$ for all $t\geq 0$ if $u\in K_{J_{m,n}}$. Define
\begin{equation}
C_k^{\left( m,n\right) }=\inf\limits_{t\geq 0}\sup\limits_{u\in \overline{%
D_{r_k,R_k}}}J_{m,n}\left( \tau \left( t,u\right) \right)  \label{eq45}
\end{equation}
with $D_{r_k,R_k}:=\left\{ u\in X_k:r_k<\left\| u\right\| _X<R_k\right\} $.

In this step we show that for $\forall t>0$,
\begin{equation}
\sup_{u\in \overline{D_{r_k,R_k}}}J_{m,n}\left( \tau \left( t,u\right)
\right) >\frac{3d_0}4  \label{eq46}
\end{equation}
for $m\in \Bbb{N}$ large enough and so $\left( J_{m,n}\right) _{\frac{3d_0}%
4}\subset X_n\backslash K_{X_n}^{2\mu }$ for $\mu >0$ sufficiently small,
where $\left( J_{m,n}\right) _\alpha :=\left\{ u\in X_n:J_{m,n}\left(
u\right) \geq \alpha \right\} $ for $\alpha \in \Bbb{R}$. Take $u_0\in
\overline{D_{r_k,R_k}}$, $\left\| u_0\right\| _X=r_k$, $u_1\in \overline{%
D_{r_k,R_k}}$, $\left\| u_1\right\| _X=R_k$, such that $\Gamma \left(
s\right) =\left( 1-s\right) u_0+su_1\subset \overline{D_{r_k,R_k}}$ for $%
\forall s\in \left( 0,1\right) $. Thereby
\begin{equation}
\tau \left( t,\Gamma \left( s\right) \right) \cap \partial B_n\left( 0,\rho
_0\right) \neq \varnothing .  \label{eq47}
\end{equation}

The boundedness of $\partial B_n\left( 0,\rho _0\right) $ in $L^\infty
\left( \Omega \right) $ follows from the fact that $X_n$ is
finite-dimensional and all norms on a finite-dimensional space are
equivalent. Hence there exists a $M>0$,
\begin{equation}
\sup\limits_{u\in \partial B_n\left( 0,\rho _0\right) }\left\| u\right\|
_{L^\infty \left( \Omega \right) }\leq M.  \label{eq48}
\end{equation}
As a consequence of $f_m\left( x,t\right) \rightarrow $\textrm{\ }$f\left(
x,t\right) $ uniformly on any compact subset of $\overline{\Omega }\times
\left[ -M,M\right] $ via $\left( \ref{eq34}\right) $, i.e., $\forall
\varepsilon >0$, $\exists m^{*}\in \Bbb{N}$, for $\forall m\geq m^{*}$, $%
\forall \sigma \in \left( 0,M\right] $, $\forall x\in \overline{\Omega }$, $%
\forall t\in \left[ -\sigma ,\sigma \right] $,
\begin{equation}
\left| f_m\left( x,t\right) -f\left( x,t\right) \right| <\varepsilon ,
\label{eq49}
\end{equation}
$\left( \ref{eq47}\right) $ and $\left( \ref{eq49}\right) $ together justify
$\left( \ref{eq46}\right) $.

At the end of Step1, we demonstrate the uniform boundedness of the minimax
value $C_k^{(m,n)}$\textbf{. }While the constant $M$ in $\left( \ref{eq48}%
\right) $ depends on the dimension $n$, the minimax value $C_k^{\left(
m,n\right) }$ is, crucially, bounded above by a constant $\beta _k$ which,
for sufficiently large $m$(depending on $n$), is ultimately independent of $%
n $ and $m$. This uniform boundedness is essential for the subsequent limit
processes in Steps 5 and 6.

To establish the uniform boundedness, we first observe that the supremum of $%
J_{m,n}$ over the initial annulus controls the minimax value:
\[
C_k^{(m,n)}\leq \sup_{u\ \in \overline{D_{r_k,R_k}}}J_{m,n}(u).
\]
Let $M^{*}:=\sup\limits_{u\in \overline{D_{r_k,R_k}}\text{ }}\left\|
u\right\| _{L^\infty }$. Since $\overline{D_{r_k,R_k}}$ is a fixed
finite-dimensional set, $M^{*}$ is finite. The uniform convergence $%
f_m\rightarrow $\textrm{\ }$f$ on any compact subset of $\overline{\Omega }%
\times \left[ -M^{*},M^{*}\right] $ implies that for any $\varepsilon >0$,
there exists a $m_{*}\in \Bbb{N}$ such that for all $m\geq m_{*}$,
\begin{equation}
\sup_{u\ \in \overline{D_{r_k,R_k}}}\left| J_{m,n}(u)-J_n(u)\right| \leq
\varepsilon .  \label{eq50}
\end{equation}
In particular, for $\varepsilon =1$ and sufficiently large $m$, we have
\[
\sup_{u\ \in \overline{D_{r_k,R_k}}}J_{m,n}(u)\leq \sup_{u\ \in \overline{%
D_{r_k,R_k}}}J_n(u)+1.
\]
Finally, since $J_n$ is simply the restriction of $J$ to $X_n$, and $%
\overline{D_{r_k,R_k}}\subset X_k$, it follows that
\[
\sup_{u\ \in \overline{D_{r_k,R_k}}}J_n(u)=\sup_{u\ \in \overline{D_{r_k,R_k}%
}}J(u)\leq \sup_{u\ \in X_k}J(u).
\]
Combining these estimates yields the desired uniform bound:
\begin{equation}
C_k^{(m,n)}\leq \sup_{u\ \in X_k}J(u)+1.  \label{eq51}
\end{equation}
The right-hand side depends only on the fixed finite-dimensional space $X_k$%
, and not on the approximation parameters $n$ or $m$. This ensures the
viability of the limit passages $m\to \infty $ and $n\to \infty $.

Step 2.\textbf{\ }In this step, we prove that for fixed $n\in \Bbb{N}$, $m>%
\widetilde{m}:=\max \left\{ m^{*},m_{*}\right\} $ and $m$ large
sufficiently, $C_k^{\left( m,n\right) }$ is a critical value of $J_{m,n}$ on
$X_n$. Assume instead that $C_k^{\left( m,n\right) }$ is a regular value.
Resting on this assumption, we claim that there exist constants $\varepsilon
_0>0$ and $t_0>0$ such that
\begin{equation}
\sup\limits_{u\in \overline{D_{r_k,R_k}}}J_{m,n}\left( \tau \left(
t_0,u\right) \right) \leq C_k^{\left( m,n\right) }-\varepsilon _0.
\label{eq52}
\end{equation}
If this were not true, then for $\forall \varepsilon >0$, $\forall t>0$,
\begin{equation}
\sup\limits_{u\in \overline{D_{r_k,R_k}}}J_{m,n}\left( \tau \left(
t,u\right) \right) >C_k^{\left( m,n\right) }-\varepsilon .  \label{eq53}
\end{equation}
Since $C_k^{\left( m,n\right) }$ is a regular value by assumption, there
exists a $\varepsilon ^{*}>0$ for which the following holds
\[
K_{J_{m,n}}\cap J_{m,n}^{-1}\left[ C_k^{\left( m,n\right) }-\varepsilon
^{*},C_k^{\left( m,n\right) }+\varepsilon ^{*}\right] =\varnothing .
\]
Hence there exists a $\gamma _0>0$ such that
\begin{equation}
\left\| J_{m,n}^{\prime }\left( u\right) \right\| _X\geq \gamma _0
\label{eq54}
\end{equation}
for $u\in J_{m,n}^{-1}\left[ C_k^{\left( m,n\right) }-\varepsilon
^{*},C_k^{\left( m,n\right) }+\varepsilon ^{*}\right] $, where $%
K_{J_{m,n}}:=\left\{ x\in X_n:J_{m,n}^{\prime }\left( x\right) =0\right\} $.

Take $\varepsilon _i>0$ and $\varepsilon _i\rightarrow 0$. Given $i\in \Bbb{N%
}$, there exists a $t_i>0$,
\[
\sup\limits_{u\in \overline{D_{r_k,R_k}}}J_{m,n}\left( \tau \left(
t_i,u\right) \right) <C_k^{\left( m,n\right) }+\varepsilon _i.
\]
Therefore, for $\forall t>0$, and $i$ large enough,
\begin{eqnarray*}
\frac{d_0}2 &<&C_k^{\left( m,n\right) }-\varepsilon _i<\sup\limits_{u\in
\overline{D_{r_k,R_k}}}J_{m,n}\left( \tau \left( t_i+t,u\right) \right) \\
&\leq &\sup\limits_{u\in \overline{D_{r_k,R_k}}}J_{m,n}\left( \tau \left(
t_i,u\right) \right) <C_k^{\left( m,n\right) }+\varepsilon _i
\end{eqnarray*}
via the combination of $\left( \ref{eq46}\right) $ and $\left( \ref{eq53}%
\right) $.

Set $J_{m,n}\left( \tau \left( t_i+t,\widetilde{u}_{t_i+t}\right) \right)
=\sup\limits_{u\in \overline{D_{r_k,R_k}}}J_{m,n}\left( \tau \left(
t_i+t,u\right) \right) $, $\widetilde{u}_{t_i+t}\in D_{r_k,R_k}$. For $%
\forall s\in \left[ t_i,t_i+t\right] $, we have
\[
\tau \left( s,\widetilde{u}_{t_i+t}\right) \in J_{m,n}^{-1}\left[
C_k^{\left( m,n\right) }-\varepsilon ^{*},C_k^{\left( m,n\right)
}+\varepsilon ^{*}\right]
\]
and consequently by $\left( \ref{eq54}\right) $,
\begin{equation}
\left\| J_{m,n}^{\prime }\left( \tau \left( s,\widetilde{u}_{t_i+t}\right)
\right) \right\| _X\geq \gamma _0.  \label{eq55}
\end{equation}
According to the definition of the flow $\left( \ref{eq44}\right) $, by
employing $\left( \ref{eq55}\right) $ we derive
\begin{eqnarray*}
&&\ J_{m,n}\left( \tau \left( t_i+t,\widetilde{u}_{t_i+t}\right) \right)
-J_{m,n}\left( \tau \left( t_i,\widetilde{u}_{t_i+t}\right) \right) \\
&=&\int_{t_i}^{t_i+t}\left\langle J_{m,n}^{\prime }\left( \tau \left( s,%
\widetilde{u}_{t_i+t}\right) \right) ,\tau _s^{\prime }\left( s,\widetilde{u}%
_{t_i+t}\right) \right\rangle ds \\
&\leq &-\int_{t_i}^{t_i+t}\left\| J_{m,n}^{\prime }\left( \tau \left( s,%
\widetilde{u}_{t_i+t}\right) \right) \right\| _X^2ds\leq -t\gamma _0^2,
\end{eqnarray*}
i.e.,
\[
J_{m,n}\left( \tau \left( t_i+t,\widetilde{u}_{t_i+t}\right) \right) \leq
J_{m,n}\left( \tau \left( t_i,\widetilde{u}_{t_i+t}\right) \right) -t\gamma
_0^2\leq C_k^{\left( m,n\right) }+\varepsilon _i-t\gamma _0^2
\]
and this yields
\[
C_k^{\left( m,n\right) }\leq \sup\limits_{u\in \overline{D_{r_k,R_k}}%
}J_{m,n}\left( \tau \left( t_i+t,u\right) \right) \leq C_k^{\left(
m,n\right) }+\varepsilon _i-t\gamma _0^2,
\]
whereupon we get a contradiction for fixed $t>0$. The claim $\left( \ref
{eq52}\right) $ is thus proved. However, $\left( \ref{eq52}\right) $
violates the definition of $C_k^{\left( m,n\right) }$. The soundness of the
argument initiated by step 2 is accordingly backed up.

Step 3. Given $t>0$, define
\[
A_t:=\left\{ v\in X_n:v=\tau \left( t,u\right) ,u\in D_{r_k,R_k}\right\}
\]
and
\[
M_t^{\left( m,n\right) }:=\left\{ \widetilde{v}\in A_t:J_{m,n}\left(
\widetilde{v}\right) =\sup\limits_{u\in \overline{D_{r_k,R_k}}}J_{m,n}\left(
\tau \left( t,u\right) \right) \right\} .
\]

Set $K_{C_k^{\left( m,n\right) }}:=\left\{ u\in K_{J_{m,n}}:J_{m,n}\left(
u\right) =C_k^{\left( m,n\right) }\right\} $. The objective of this step is
to establish that there exist $\beta _j^{*}>0$, $\beta _j^{*}\rightarrow 0$,
and $t_j>0$, such that
\begin{equation}
\text{dist}\left( M_{t_j,\beta _j^{*}}^{\left( m,n\right) },K_{C_k^{\left(
m,n\right) }}\right) \rightarrow 0  \label{eq56}
\end{equation}
as $j\rightarrow \infty $, where
\[
M_{t,\beta }^{\left( m,n\right) }:=\left\{ \widetilde{v}\in M_t^{\left(
m,n\right) }:J_{m,n}\left( \widetilde{v}\right) <C_k^{\left( m,n\right)
}+\beta \right\}
\]
for $t>0$, $\beta >0$.

Notice that for fixed $t\geq 0$, $\tau :\left( t,\cdot \right) \rightarrow
A_t$ is a $C^2$ diffeomorphism for $\forall u\in D_{r_k,R_k}$, therefore for
$v\in A_t\subset X_n$, $\exists |$ $u\in D_{r_k,R_k}$, $v=\tau \left(
t,u\right) $, and we set $u:=\tau _t^{-1}\left( v\right) $, $\tau _t\left(
\cdot \right) :=\tau \left( t,\cdot \right) $ for $u\in X_n$. We will verify
that $\forall \beta >0$, $\forall \delta >0$, $\forall t>0$, $\exists
t^{*}>t $, $\exists v_{t^{*}}\in M_{t^{*},\beta }^{\left( m,n\right) }$, $%
\exists \widetilde{t}^{*}\in \left[ t^{*}-t,t^{*}\right] $, such that
\begin{equation}
\left\| J_{m,n}^{\prime }\left( \tau \left( \widetilde{t}^{*},\tau
_{t^{*}}^{-1}\left( v_{t^{*}}\right) \right) \right) \right\| _X<\delta .
\label{eq57}
\end{equation}

By way of negation, $\exists \beta _0>0$, $\exists \delta _0>0$, $\exists
t_0>0$, $\forall t>t_0$, $\forall v_t\in M_{t,\beta _0}^{\left( m,n\right) }$%
, $\forall \widetilde{t}\in \left[ t-t_0,t\right] $,
\begin{equation}
\left\| J_{m,n}^{\prime }\left( \tau \left( \widetilde{t},\tau _t^{-1}\left(
v_t\right) \right) \right) \right\| _X\geq \delta _0.  \label{eq58}
\end{equation}
Consequently, $\left( \ref{eq44}\right) $ together with $\left( \ref{eq46}%
\right) $ gives rise to
\begin{eqnarray}
&&\ \ \ J_{m,n}\left( \tau \left( t,\tau _t^{-1}\left( v_t\right) \right)
\right) -J_{m,n}\left( \tau \left( t-t_0,\tau _t^{-1}\left( v_t\right)
\right) \right)  \nonumber \\
\ &=&\int_{t-t_0}^t\left\langle J_{m,n}^{\prime }\left( \tau \left( s,\tau
_t^{-1}\left( v_t\right) \right) \right) ,\tau _s^{\prime }\left( s,\tau
_t^{-1}\left( v_t\right) \right) \right\rangle _{X_n}ds  \nonumber \\
\ &\leq &-\int_{t-t_0}^t\left\| J_{m,n}^{\prime }\left( \tau \left( s,\tau
_t^{-1}\left( v_t\right) \right) \right) \right\| _X^2ds\leq -\delta _0^2t_0,
\label{eq59}
\end{eqnarray}
i.e.,
\begin{equation}
J_{m,n}\left( v_t\right) =J_{m,n}\left( \tau \left( t,\tau _t^{-1}\left(
v_t\right) \right) \right) \leq J_{m,n}\left( \tau \left( t-t_0,\tau
_t^{-1}\left( v_t\right) \right) \right) -\delta _0^2t_0.  \label{eq60}
\end{equation}
We can therefore find sequences $\widetilde{t}_j\rightarrow \infty $ and $%
\beta _j\rightarrow 0$, with $v_{\widetilde{t}_j}\in M_{\widetilde{t}%
_j,\beta _j}^{\left( m,n\right) }$. Setting $t_j:=\widetilde{t}_j+t_0$, we
obtain a contradiction
\begin{eqnarray}
C_k^{\left( m,n\right) } &\leq &J_{m,n}\left( v_{t_j}\right) \leq
J_{m,n}\left( \tau \left( \widetilde{t}_j,\tau _{t_j}^{-1}\left(
v_{t_j}\right) \right) \right) -\delta _0^2t_0  \nonumber \\
&\leq &C_k^{\left( m,n\right) }+\beta _j-\delta _0^2t_0\rightarrow
C_k^{\left( m,n\right) }-\delta _0^2t_0.  \label{eq61}
\end{eqnarray}
The assertion $\left( \ref{eq57}\right) $ follows.

Let $\widetilde{\beta }_i\rightarrow 0$, $\delta _i\rightarrow 0$, $%
t_i\rightarrow 0$. According to $\left( \ref{eq57}\right) $, for fixed $i\in
\Bbb{N}$, there exists a $t_i^{*}>t_i$, $\exists v_{t_i^{*}}\in M_{t_i^{*},%
\widetilde{\beta }_i}^{\left( m,n\right) }$, $\exists \widetilde{t}_i^{*}\in
\left[ t_i^{*}-t_i,t_i^{*}\right] $,
\begin{equation}
\left\| J_{m,n}^{\prime }\left( \tau \left( \widetilde{t}_i^{*},\tau
_{t_i^{*}}^{-1}\left( v_{t_i^{*}}\right) \right) \right) \right\| _X<\delta
_i.  \label{eq62}
\end{equation}
Set $u_i=\tau _{t_i^{*}}^{-1}\left( v_{t_i^{*}}\right) $. We claim that $%
\exists M>0$, for $i\in \Bbb{N}$,
\begin{equation}
\sup\limits_{\widetilde{t}\in \left[ t_i^{*}-t_i,t_i^{*}\right] }\left\|
\tau \left( \widetilde{t},u_i\right) \right\| _X\leq M.  \label{eq63}
\end{equation}
Otherwise, there exists $\zeta _{i_j}\in \left[
t_{i_j}^{*}-t_{i_j},t_{i_j}^{*}\right] $, $\left\| \tau \left( \zeta
_{i_j},u_{i_j}\right) \right\| _X\rightarrow \infty $ as $j\rightarrow
\infty $. Set $\sigma _{i_j}=\frac{\tau \left( \zeta _{i_j},u_{i_j}\right) }{%
\left\| \tau \left( \zeta _{i_j},u_{i_j}\right) \right\| _X}$, $\tau \left(
\zeta _{i_j},u_{i_j}\right) =\gamma _{i_j}\sigma _{i_j}$, $\gamma _{i_j}>0$.
The unboundedness of $\left\{ \left\| \tau \left( \zeta
_{i_j},u_{i_j}\right) \right\| _X\right\} _{i=1}^\infty $ gives rise to an
absurd deduction
\begin{eqnarray}
C_k^{\left( m,n\right) } &\leq &J_{m,n}\left( \tau \left( \zeta
_{i_j},u_{i_j}\right) \right) =\frac 12\left\| \tau \left( \zeta
_{i_j},u_{i_j}\right) \right\| _X^2-\int_\Omega F_m\left( x,\tau \left(
\zeta _{i_j},u_{i_j}\right) \right)  \nonumber \\
\ &\leq &\frac 12\gamma _{i_j}^2-\widetilde{C}_1\gamma _{i_j}^{\widetilde{%
\theta }}\int_\Omega \left| \sigma _{i_j}\right| ^{\widetilde{\theta }}+%
\widetilde{C}_2\leq \frac 12\gamma _{i_j}^2-\widetilde{C}_3\gamma _{i_j}^{%
\widetilde{\theta }}+\widetilde{C}_2\rightarrow -\infty ,  \label{eq64}
\end{eqnarray}
where $\widetilde{C}_1$,$\widetilde{C}_2$,$\widetilde{C}_3>0$ denote
constants independent of $m$. We conclude the claim $\left( \ref{eq63}%
\right) $ as desired. By the integral mean value theorem,
\begin{eqnarray}
&&\ J_{m,n}^{\prime }\left( \tau \left( t_i^{*},u_i\right) \right)
-J_{m,n}^{\prime }\left( \tau \left( \widetilde{t}_i^{*},u_i\right) \right)
\ \ \ \ \   \nonumber \\
&=&\int_0^1\frac d{d\xi }J_{m,n}^{\prime }\left( y_{\xi ,i}\right) d\xi
\nonumber \\
\ &=&\left( t_i^{*}-\widetilde{t}_i^{*}\right) \int_0^1J_{m,n}^{^{\prime
\prime }}\left( y_{\xi ,i}\right) \tau _t^{\prime }\left( \xi t_i^{*}+\left(
1-\xi \right) \widetilde{t}_i^{*},u_i\right) d\xi  \nonumber \\
\ &=&-\left( t_i^{*}-\widetilde{t}_i^{*}\right) \cdot
\int_0^1J_{m,n}^{^{\prime \prime }}\left( y_{\xi ,i}\right) \left( 1-l\left(
y_{\xi ,i}\right) \right) J_{m,n}^{\prime }\left( y_{\xi ,i}\right) d\xi ,
\label{eq65}
\end{eqnarray}
where $y_{\xi ,i}:=\tau \left( \xi t_i^{*}+\left( 1-\xi \right) \widetilde{t}%
_i^{*},u_i\right) $. Set $\Theta _{m,n}:=\left\{ u\in X_n:u=y_{\xi ,i}\text{%
, }\xi \in \left[ 0,1\right] \text{, }i\in \Bbb{N}\right\} $. Notice that
the closure of $\Theta _{m,n}$, denoted by $\overline{\Theta _{m,n}}$, is a
compact set in $X_n$ since $\left( \ref{eq63}\right) $ yields
\begin{equation}
\sup\limits_{u\in \overline{\Theta _{m,n}}}\left\| u\right\| _X\leq M.
\label{eq66}
\end{equation}
Hence, there exist constants $M_0$,$M_1>0$ such that
\[
\sup\limits_{u\in \overline{\Theta _{m,n}}}\left\| J_{m,n}^{^{\prime \prime
}}\left( u\right) \right\| _{L\left( X\right) }\leq M_0;\text{ }%
\sup\limits_{u\in \overline{\Theta _{m,n}}}\left\| J_{m,n}^{\prime }\left(
u\right) \right\| _X\leq M_1.
\]
Consequently, $\left( \ref{eq65}\right) $ together with $\left( \ref{eq66}%
\right) $ yields
\begin{eqnarray}
&&\ \ \ \ \ \ \ \left\| J_{m,n}^{\prime }\left( \tau \left(
t_i^{*},u_i\right) \right) -J_{m,n}^{\prime }\left( \tau \left( \widetilde{t}%
_i^{*},u_i\right) \right) \right\| _X  \nonumber \\
\ &\leq &\left( t_i^{*}-\widetilde{t}_i^{*}\right) \int_0^1\left\|
J_{m,n}^{^{\prime \prime }}\left( y_{\xi ,i}\right) J_{m,n}^{\prime }\left(
y_{\xi ,i}\right) \right\| _Xd\xi \leq \left( t_i^{*}-\widetilde{t}%
_i^{*}\right) M_0M_1\rightarrow 0.  \label{eq67}
\end{eqnarray}
The anticipation $\left( \ref{eq56}\right) $ is achieved.

Step 4. This step demonstrates that for $k\in \Bbb{N}$ and $k<n$, $J_{m,n}$
has a critical point $u$ on $X_n$ with $M\left( J_{m,n},u\right) \geq k$.
Employing Lemma \ref{lem2.1}, it is clear that the map $\tau \left( t,\cdot
\right) $ is a $C^2$ diffeomorphism onto its image, hence $A_t$ is, by
definition, a $k$-dimensional $C^2$ Riemannian submanifold of $X_n$.

Combining $\left( \ref{eq62}\right) $ and $\left( \ref{eq67}\right) $,
\begin{equation}
\left\| J_{m,n}^{\prime }\left( v_{t_i^{*}}\right) \right\| _X=\left\|
J_{m,n}^{\prime }\left( \tau \left( t_i^{*},u_i\right) \right) \right\|
_X\leq \left\| J_{m,n}^{\prime }\left( \tau \left( \widetilde{t}%
_i^{*},u_i\right) \right) \right\| _X+o_i\rightarrow 0,  \label{eq68}
\end{equation}
where
\begin{equation}
o_i=\ \left\| J_{m,n}^{\prime }\left( \tau \left( t_i^{*},u_i\right) \right)
-J_{m,n}^{\prime }\left( \tau \left( \widetilde{t}_i^{*},u_i\right) \right)
\right\| _{_X}.  \label{eq69}
\end{equation}
Due to
\begin{equation}
J_{m,n}\left( \tau \left( t_i^{*},u_i\right) \right) =\sup_{u\in
D_{r_k,R_k}}J_{m,n}\left( \tau \left( t_i^{*},u\right) \right) ,
\label{eq70}
\end{equation}
we get
\begin{equation}
\left\langle \left[ J_{m,n}\left( \tau \left( t_i^{*},u\right) \right)
\right] _{u=u_i}^{\prime },\varphi \right\rangle _{X_k}=\left\langle \tau
_u^{\prime }\left( t_i^{*},u_i\right) ^{*}J_{m,n}^{\prime }\left( \tau
\left( t_i^{*},u_i\right) \right) ,\varphi \right\rangle _{X_k}=0
\label{eq71}
\end{equation}
for $\forall \varphi \in X_k$, where $\tau _u^{\prime }\left(
t_i^{*},u_i\right) ^{*}$ represents the conjugate transpose operator of $%
\tau _u^{\prime }\left( t_i^{*},u_i\right) $. Let $T_{\tau \left( t,u\right)
}A_t$ be the tangent space of $A_t$ at $\tau \left( t,u\right) $ for $t\geq
0 $ and $u\in D_{r_k,R_k}$. Observe that
\begin{equation}
\left( \tau _u^{\prime }\left( t_i^{*},u_i\right) ^{*}\right) ^{-1}=\left(
\tau _u^{\prime }\left( t_i^{*},u_i\right) ^{-1}\right) ^{*},  \label{eq72}
\end{equation}
and also
\begin{eqnarray}
\tau _u^{\prime }\left( t_i^{*},u_i\right) &\in &L\left( X_k,T_{\tau \left(
t_i^{*},u_i\right) }A_{t_i^{*}}\right) ,  \label{eq73} \\
\tau _u^{\prime \prime }\left( t_i^{*},u_i\right) &\in &L\left( X_k,L\left(
X_k,T_{\tau \left( t_i^{*},u_i\right) }A_{t_i^{*}}\right) \right) ,
\label{eq74}
\end{eqnarray}
using $\left( \ref{eq71}\right) $ and $\left( \ref{eq72}\right) $ we
conclude that for $\forall \varphi $,$\psi \in X_k$,
\begin{eqnarray}
&&\left\langle \left[ J_{m,n}\left( \tau \left( t_i^{*},u\right) \right)
\right] _{u=u_i}^{\prime \prime }\varphi ,\psi \right\rangle _{X_n}
\nonumber \\
&=&\left\langle \left[ \tau _u^{\prime }\left( t_i^{*},u\right)
^{*}J_{m,n}^{\prime }\left( \tau \left( t_i^{*},u\right) \right) \right]
_{u=u_i}^{\prime }\varphi ,\psi \right\rangle _{X_n}  \nonumber \\
\ &=&\left\langle \left[ \tau _u^{\prime \prime }\left( t_i^{*},u_i\right)
\varphi \right] ^{*}\left( \tau _u^{\prime }\left( t_i^{*},u_i\right)
^{*}\right) ^{-1}\tau _u^{\prime }\left( t_i^{*},u_i\right)
^{*}J_{m,n}^{\prime }\left( \tau \left( t_i^{*},u_i\right) \right) ,\psi
\right\rangle _{X_n}  \nonumber \\
&&\ +\left\langle \tau _u^{\prime }\left( t_i^{*},u_i\right)
^{*}J_{m,n}^{\prime \prime }\left( \tau \left( t_i^{*},u_i\right) \right)
\tau _u^{\prime }\left( t_i^{*},u_i\right) \varphi ,\psi \right\rangle _{X_n}
\nonumber \\
&=&\left\langle \left( \tau _u^{\prime }\left( t_i^{*},u_i\right)
^{-1}\right) ^{*}\tau _u^{\prime }\left( t_i^{*},u_i\right)
^{*}J_{m,n}^{\prime }\left( \tau \left( t_i^{*},u_i\right) \right) ,\left[
\tau _u^{\prime \prime }\left( t_i^{*},u_i\right) \varphi \right] \psi
\right\rangle _{X_n}  \nonumber \\
&&\ +\left\langle J_{m,n}^{\prime \prime }\left( \tau \left(
t_i^{*},u_i\right) \right) \tau _u^{\prime }\left( t_i^{*},u_i\right)
\varphi ,\tau _u^{\prime }\left( t_i^{*},u_i\right) \psi \right\rangle _{X_n}
\nonumber \\
&=&\left\langle \tau _u^{\prime }\left( t_i^{*},u_i\right)
^{*}J_{m,n}^{\prime }\left( \tau \left( t_i^{*},u_i\right) \right) ,\tau
_u^{\prime }\left( t_i^{*},u_i\right) ^{-1}\left[ \tau _u^{\prime \prime
}\left( t_i^{*},u_i\right) \varphi \right] \psi \right\rangle _{X_n}
\nonumber \\
&&+\left\langle J_{m,n}^{\prime \prime }\left( \tau \left(
t_i^{*},u_i\right) \right) \tau _u^{\prime }\left( t_i^{*},u_i\right)
\varphi ,\tau _u^{\prime }\left( t_i^{*},u_i\right) \psi \right\rangle _{X_n}
\nonumber \\
\ &=&\left\langle J_{m,n}^{\prime \prime }\left( \tau \left(
t_i^{*},u_i\right) \right) \tau _u^{\prime }\left( t_i^{*},u_i\right)
\varphi ,\tau _u^{\prime }\left( t_i^{*},u_i\right) \psi \right\rangle
_{X_n}.  \label{eq75}
\end{eqnarray}

By $\left( \ref{eq70}\right) $ we have
\begin{equation}
\left\langle \left[ J_{m,n}\left( \tau \left( t_i^{*},u\right) \right)
\right] _{u=u_i}^{\prime \prime }\varphi ,\varphi \right\rangle _{X_n}\leq 0.
\label{eq76}
\end{equation}
Hence by employing $\left( \ref{eq75}\right) $ and $\left( \ref{eq76}\right)
$,
\begin{eqnarray}
&&\ \sup_{\varphi \in X_k,\left\| \varphi \right\| _X=1}\left\langle
J_{m,n}^{\prime \prime }\left( \tau \left( t_i^{*},u_i\right) \right) \tau
_u^{\prime }\left( t_i^{*},u_i\right) \varphi ,\tau _u^{\prime }\left(
t_i^{*},u_i\right) \varphi \right\rangle _{X_n}  \nonumber \\
\ &=&\sup_{\varphi \in X_k,\left\| \varphi \right\| _X=1}\left\langle \left[
\tau ^{\prime }\left( t_i^{*},u\right) ^{*}J_{m,n}^{\prime }\left( \tau
\left( t_i^{*},u\right) \right) \right] _{u=u_i}^{\prime }\varphi ,\varphi
\right\rangle _{X_n}  \nonumber \\
\ &=&\sup_{\varphi \in X_k,\left\| \varphi \right\| _X=1}\left\langle \left[
J_{m,n}\left( \tau \left( t_i^{*},u\right) \right) \right] _{u=u_i}^{\prime
\prime }\varphi ,\varphi \right\rangle _{X_n}\leq 0.  \label{eq77}
\end{eqnarray}
Consequently $\left( \ref{eq77}\right) $ derives
\begin{equation}
\sup_{\left\| \psi \right\| _X=1,\psi \in T_{\tau \left( t_i^{*},u_i\right)
}A_{t_i^{*}}}\left\langle J_{m,n}^{\prime \prime }\left( \tau \left(
t_i^{*},u_i\right) \right) \psi ,\psi \right\rangle _{X_n}\leq 0.
\label{eq78}
\end{equation}

It is noteworthy that $\left\{ \tau \left( t_i^{*},u_i\right) \right\}
_{i=1}^\infty $ is a bounded (PS) sequence for $J_{m,n}$ in $X_n$ by $\left(
\ref{eq68}\right) $. Therefore, by $\left( \ref{eq51}\right) $, a
subsequence of $\tau \left( t_i^{*},u_i\right) $ converges to a critical
point $u_{m,n,k}$ of $J_{m,n}$ in $X_n$,
\begin{equation}
J_{m,n}\left( u_{m,n,k}\right) =C_k^{\left( m,n\right) }\leq \sup_{u\in
X_k}J\left( u\right) +1  \label{eq79}
\end{equation}
for all sufficiently large $m$.

The remainder of the proof of Step 4 aims to argue the claim $M\left(
J_{m,n},u_{m,n,k}\right) \geq k$. By way of negation, suppose that $\dim
\left( X_n^{+}\left( u_{m,n,k}\right) \right) ^{\perp _n}<k$, there exist
\begin{equation}
\psi _{t_i^{*}}\in \left( X_n^{-}\left( \tau \left( t_i^{*},u_i\right)
\right) \oplus X_n^0\left( \tau \left( t_i^{*},u_i\right) \right) \right)
\cap X_n^{+}\left( u_{m,n,k}\right) ,  \label{eq80}
\end{equation}
$\left\| \psi _{t_i^{*}}\right\| _X=1$, as $t_i^{*}\rightarrow \infty $,
where $X_n^{+}\left( u_{m,n,k}\right) $ is the positive eigenspace of $%
J_{m,n}^{\prime \prime }\left( u_{m,n,k}\right) $ on $X_n$ and $\left(
X_n^{+}\left( u_{m,n,k}\right) \right) ^{\perp _n}$ stands for the
orthogonal complement in $X_n$ of $X_n^{+}\left( u_{m,n,k}\right) $, and
additionally $X_n^{-}\left( \tau \left( t_i^{*},u_i\right) \right) $,$%
X_n^0\left( \tau \left( t_i^{*},u_i\right) \right) $ are the negative and
kernel subspaces of $J_{m,n}^{\prime \prime }\left( \tau \left(
t_i^{*},u_i\right) \right) $ respectively. Based on the definition of $%
X_n^{+}\left( u_{m,n,k}\right) $, $\exists \rho >0$, such that
\[
\inf\limits_{v\in X_n^{+}\left( u_{m,n,k}\right) \cap S^{n-1}}\left\langle
J_{m,n}^{\prime \prime }\left( u_{m,n,k}\right) v,v\right\rangle _{X_n}>\rho
,
\]
where $S^{n-1}=\left\{ v\in X_n:\left\| v\right\| _X=1\right\} $. Notice
that $J_{m,n}^{\prime \prime }\left( \tau \left( t_i^{*},u_i\right) \right)
\rightarrow J_{m,n}^{\prime \prime }\left( u_{m,n,k}\right) $ as $%
t_i^{*}\rightarrow \infty $. Pick $t_i^{*}$ sufficiently large such that
\[
\left| \left\langle \left( J_{m,n}^{\prime \prime }\left( \tau \left(
t_i^{*},u_i\right) \right) -J_{m,n}^{\prime \prime }\left( u_{m,n,k}\right)
\right) \psi _{t_i^{*}},\psi _{t_i^{*}}\right\rangle _{X_n}\right| <\frac
\rho 2,
\]
and hence we get an absurd inequality
\begin{equation}
0\geq \left\langle J_{m,n}^{\prime \prime }\left( \tau \left(
t_i^{*},u_i\right) \right) \psi _{t_i^{*}},\psi _{t_i^{*}}\right\rangle
_{X_n}>\frac \rho 2,  \label{eq81}
\end{equation}
subverting the contradiction hypothesis $\dim \left( X_n^{+}\left(
u_{m,n,k}\right) \right) ^{\perp _n}<k$. We consequently confirm the
assertion $M\left( J_{m,n},u_{m,n,k}\right) \geq k$.

Step 5.\textbf{\ }Let $J_n^{\prime }\left( u\right) :=P_nJ^{\prime }\left(
u\right) $ for $u\in X$ and $K_{J_n}:=\left\{ x\in X_n:J_n^{\prime }\left(
x\right) =0\right\} $. In this step we establish the boundedness of $\left\{
\left\| u_{m,n,k}\right\| _X\right\} _{m=1}^\infty $ and then prove that $%
\left\{ u_{m,n,k}\right\} _{m=1}^\infty $ converges to a critical point of $%
J_n$ in $X_n$.

The following definition and lemma are adapted to our needs:

\begin{definition}
\label{def4.3} For fixed $n\in \Bbb{N}$, we say that $\left\{
J_{m,n}\right\} _{m=1}^\infty $ satisfies the (PS)$^{\text{**}}$ condition
on $X_n$ if any sequence $\left\{ u_m\right\} _{m=1}^\infty \subset X_n$
with $\left\{ J_{m,n}\left( u_m\right) \right\} _{m=1}^\infty $ bounded from
above and $\left\| J_{m,n}^{\prime }\left( u_m\right) \right\| _X\rightarrow
0$ contains a subsequence convergent to a critical point of $J_n$ on $X_n$.
\end{definition}

\begin{lemma}
\label{lem4.4} Suppose that $\left( f_1\right) $ and $\left( f_2\right) $
hold, then for fixed $n\in \Bbb{N}$, $\left\{ J_{m,n}\right\} _{m=1}^\infty $
satisfies the (PS)$^{\text{**}}$ condition on $X_n$.
\end{lemma}

\textit{Proof}. Assumption $\left( f_1\right) $  holds for $f_m$ for all
sufficiently large $m\in \Bbb{N}$. Moreover, Lemma \ref{lem3.2} implies that
assumption $\left( f_2\right) $ is also satisfied by $f_m$ for $\widetilde{%
\theta }\in \left( 2,\theta \right) $ and $m\in \Bbb{N}$ large. Given any
sequence $\left\{ u_m\right\} _{m=1}^\infty \subset X_n$ such that $%
\sup\limits_{m\in \Bbb{N}}J_{m,n}\left( u_m\right) <\infty $ and $\left\|
J_{m,n}^{\prime }\left( u_m\right) \right\| _X\rightarrow 0$ as $%
m\rightarrow \infty $, we have
\begin{equation}
\left\| u_m\right\| _X^2=\int_\Omega f_m\left( x,u_m\right) u_mdx+o\left(
1\right) \left\| u_m\right\| _X.  \label{eq82}
\end{equation}
Argue indirectly and suppose that $\left\| u_m\right\| _X\rightarrow \infty $
as $m\rightarrow \infty $, therefore $\left( \ref{eq82}\right) $ and $\left(
f_2\right) $ together lead to the expected contradiction
\begin{eqnarray}
\ \infty  &>&\sup\limits_{m\in \Bbb{N}}J_{m,n}\left( u_m\right) \geq \frac
12\left\| u_m\right\| _X^2-\frac 1{\widetilde{\theta }}\int_\Omega f_m\left(
x,u_m\right) u_mdx  \nonumber \\
&&\ \ \ \ \ \ \ \ \ \ +\int_{\left| u_m\right| \leq M}\left( \frac 1{%
\widetilde{\theta }}f_m\left( x,u_m\right) u_m-F_m\left( x,u_m\right)
\right) dx  \nonumber \\
\  &=&\left( \frac 12-\frac 1{\widetilde{\theta }}\right) \left\|
u_m\right\| _X^2+o\left( 1\right) \frac 1{\widetilde{\theta }}\left\|
u_m\right\| _X  \nonumber \\
&&\ \ \ \ \ \ \ \ \ \ +\int_{\left| u_m\right| \leq M}\left( \frac 1{%
\widetilde{\theta }}f_m\left( x,u_m\right) u_m-F_m\left( x,u_m\right)
\right) dx  \nonumber \\
\  &\rightarrow &\infty   \label{eq83}
\end{eqnarray}
by sending $m$ to infinity, which is plainly attributable to the initial
assumption and thus yield the boundedness of $\left\{ \left\| u_m\right\|
_X\right\} _{m=1}^{+\infty }$. Let $u_m\rightarrow u_0$ in $X_n$. We remain
to verify $u_0\in K_{J_n}$. Observe that for $\forall v\in X_n$,
\begin{equation}
\left\langle \left( J_{m,n}^{\prime }\left( u_m\right) -J_n^{\prime }\left(
u_0\right) \right) ,v\right\rangle _X=\left\langle u_m-u_0,v\right\rangle
_X-\int_\Omega \left( f_m\left( x,u_m\right) -f\left( x,u_0\right) \right)
vdx.\   \label{eq84}
\end{equation}

Given $v\in X$ and pick $r_1\in \left( \frac{2N}{N+2},\frac{2^{*}}{p-1}%
\right) $, $r_2\in \left( 1,2^{*}\right) $, $\frac 1{r_1}+\frac 1{r_2}=1$.
Invoking the integral mean value theorem, for fixed $m\in \Bbb{N}$, there
exists a constant $\tau _m\in \left( -\frac 1m,\frac 1m\right) $ such that
\begin{equation}
\int_{-\frac 1m}^{\frac 1m}\rho _m\left( \tau \right) \left[ \left( f\left(
x,u_m-\tau \right) -f\left( x,u_0\right) \right) \right] d\tau =f\left(
x,u_m-\tau _m\right) -f\left( x,u_0\right) .  \label{eq85}
\end{equation}
Hence, employing Proposition \ref{prop4.2} we get
\begin{eqnarray}
&&\ \ \ \ \ \ \ \ \ \ \ \ \ \ \ \ \ \left| \int_\Omega \left( f_m\left(
x,u_m\right) -f\left( x,u_0\right) \right) vdx\right|  \nonumber \\
\ &=&\left| \int_\Omega \int_{-\frac 1m}^{\frac 1m}\rho _m\left( \tau
\right) \left[ \left( f\left( x,u_m-\tau \right) -f\left( x,u_0\right)
\right) \right] d\tau \cdot vdx\right|  \nonumber \\
\ &=&\left| \int_\Omega \left[ \left( f\left( x,u_m-\tau _m\right) -f\left(
x,u_0\right) \right) \right] vdx\right|  \nonumber \\
\ &\leq &\left( \int_\Omega \left| f\left( x,u_m-\tau _m\right) -f\left(
x,u_0\right) \right| ^{r_1}dx\right) ^{\frac 1{r_1}}\left( \int_\Omega
\left| v\right| ^{r_2}dx\right) ^{\frac 1{r_2}}  \nonumber \\
\ &\leq &C\left( \int_\Omega \left| f\left( x,u_m-\tau _m\right) -f\left(
x,u_0\right) \right| ^{r_1}dx\right) ^{\frac 1{r_1}}\cdot \left\| v\right\|
_X\rightarrow 0  \label{eq86}
\end{eqnarray}
by virtue of
\begin{equation}
\left\| u_m-\tau _m-u_0\right\| _{L^{\left( p-1\right) r_1}}\leq \left\|
u_m-u_0\right\| _{L^{\left( p-1\right) r_1}}+\left\| \tau _m\right\|
_{L^{\left( p-1\right) r_1}}\rightarrow 0  \label{eq87}
\end{equation}
\textrm{\ }as\textrm{\ }$m\rightarrow \infty $.

By inserting $\left( \ref{eq86}\right) $ into $\left( \ref{eq84}\right) $ we
derive $J_n^{\prime }\left( u_0\right) =0$ and thus have corroborated the
(PS)$^{\text{**}}$ condition with respect to $\left\{ J_{m,n}\right\}
_{m=1}^\infty $. The proof of Lemma \ref{lem4.4} is complete.\qed\vskip 5pt

We now proceed with the proof of Step 5. An argument analogous to $\left(
\ref{eq83}\right) $ shows the boundedness of $\left\{ \left\|
u_{m,n,k}\right\| _X\right\} _{m=1}^\infty $ by $\left( \ref{eq79}\right) $.
Likewise $\lim\limits_{m\rightarrow \infty }\left\| u_{m,n,k}\right\|
_{L^\infty \left( \Omega \right) }<\infty $. Let $M_{n,k}:=\sup\limits_{m\in
\Bbb{N}}\left\| u_{m,n,k}\right\| _{L^\infty \left( \Omega \right) }$. The
uniform convergence $f_m\rightarrow $\textrm{\ }$f$ on the compact set $%
\overline{\Omega }\times \left[ -M_{n,k},M_{n,k}\right] $ implies that for
any $\varepsilon >0$, there exists a $\widetilde{m}\in \Bbb{N}$ such that
for $\forall m\geq \widetilde{m}$, $\forall x\in \overline{\Omega }$, $%
\forall t\in \left[ -M_{n,k},M_{n,k}\right] $,
\[
\left| f_m\left( x,t\right) -f\left( x,t\right) \right| \leq \varepsilon .
\]
Consequently
\begin{eqnarray*}
&&\ \left| J_{m,n}(u_{m,n,k})-J_n(u_{m,n,k})\right| \\
&\leq &\int_\Omega \int_0^{M_{n,k}}\left| f_m\left( x,s\right) -f\left(
x,s\right) \right| dsdx\leq \varepsilon M_{n,k}\left| \Omega \right| .
\end{eqnarray*}
In particular, for $\varepsilon =\frac 1{M_{n,k}\left| \Omega \right| }$ and
sufficiently large $m$, by $\left( \ref{eq79}\right) $ we get
\[
J_n\left( u_{m,n,k}\right) \leq J_{m,n}\left( u_{m,n,k}\right) +1\leq
\sup_{u\in X_k}J\left( u\right) +2
\]
for $m\in \Bbb{N}$ large. Therefore $\left\{ u_{m,n,k}\right\} _{m=1}^\infty
$ converges to a critical point $u_{n,k}$ of $J_n$ in $X_n$ with
\begin{equation}
J_n\left( u_{n,k}\right) \leq \sup_{u\in X_k}J\left( u\right) +2
\label{eq88}
\end{equation}
and
\begin{equation}
M\left( J_n,u_{n,k}\right) \geq k.  \label{eq89}
\end{equation}

Step 6.\textbf{\ }Since $J$ satisfies the (PS)$^{*}$ condition with respect
to $\left\{ X_n\right\} _{n=1}^\infty $, combining $\left( \ref{eq88}\right)
$, $\left( \ref{eq89}\right) $ and $J_n^{\prime }\left( u_{n,k}\right) =0$,
we yield $u_{n,k}\rightarrow u_k\in K_J:=\left\{ u\in X:J^{\prime }\left(
u\right) =0\right\} $ with $M\left( J,u_k\right) \geq k$.

Define $k_0:=M\left( J,u_k\right) $ and take $k=k_i:=k_0+i$, $i\in \Bbb{N}$,
and set $X_{k_i}:=$span$\left\{ \varphi _1,\cdots ,\varphi _{k_i}\right\} $.
Repeat the procedure of the previous argument, we ultimately get a sequence
of solutions $\left\{ u_{k_i}\right\} _{i=1}^{+\infty }$ of the equation $%
\left( \ref{eq1}\right) $ with
\begin{equation}
M\left( J\text{, }u_{k_i}\right) \geq k_i.  \label{eq90}
\end{equation}

Thus far we have validated the promised proof of the first part of Theorem
\ref{thm1.1}. Now we are left to verify $J\left( u_{k_i}\right) \rightarrow
\infty $ as $i\rightarrow \infty $. Suppose to the contrary that there
exists a constant $\beta ^{*}>0$ such that $J\left( u_{k_i}\right) \leq
\beta ^{*}$. A standard argument derives the boundedness of $\left\{ \left\|
u_{k_i}\right\| _X\right\} _{i=1}^{+\infty }$ and assume $u_{k_i}\rightarrow
\widetilde{u}\in K_J$. As a consequence of the compactness of the operator $%
\left( -\Delta \right) ^{-1}f^{\prime }\left( u\right) :X\rightarrow X$ for
any $u\in X$, the sum of the dimensions of the negative and zero subspaces
of the Hessian $J^{\prime \prime }\left( u\right) =I-\left( -\Delta \right)
^{-1}f^{\prime }\left( u\right) $ is finite. We thereby get a contradiction
\begin{equation}
\infty \leftarrow M\left( J,u_{k_i}\right) \leq M\left( J,\widetilde{u}%
\right) <\infty  \label{eq91}
\end{equation}
via $\left( \ref{eq90}\right) $. The assertion $J\left( u_{k_i}\right)
\rightarrow \infty $ follows.\ This substantiates the statement $\left(
ii\right) $ of Theorem \ref{thm1.1}, thus completing the proof.

Indeed, we can extract an abstract theorem from the proof of Theorem \ref
{thm1.1}, which can be stated as follows:

\begin{theorem}
\label{thm4.5} Let $E$ be a Hilbert space endowed with the norm $\left\|
\cdot \right\| _E$, $\dim E=n<\infty $. Suppose that $\Phi \in C^3\left( E,%
\Bbb{R}\right) $ with $\Phi \left( 0\right) =0$ satisfies the (PS)
condition. Assume further that:

$\left\langle 1\right\rangle $ There exists an $r>0$, s.t. $\Phi \left(
x\right) >0$ for any $x\in \overline{B_E\left( 0,r\right) }\backslash
\left\{ 0\right\} $, $B_E\left( 0,r\right) :=\left\{ x\in E:\left\|
x\right\| _E<r\right\} $;

$\left\langle 2\right\rangle $ There exist $\beta >0$ and $\widetilde{R}>0$
such that $\Phi \left( x\right) \leq -\beta $ for all $x\in E$ with $\left\|
x\right\| _E=\widetilde{R}$;

$\left\langle 3\right\rangle $ The set $\left\{ x\in E:\Phi \left( x\right)
\geq \sigma \right\} $ is bounded for any $\sigma \in \Bbb{R}$.

Let $\widetilde{E}$ be a Hilbert subspace of $E$ with $\dim \widetilde{E}%
=l<n $, and suppose $\sup\limits_{x\in \widetilde{E}}\Phi \left( x\right)
<\infty $. Choose constants $0<\widetilde{r}<r$, and define the annulus:
\[
D_{\widetilde{r},\widetilde{R}}:=\left\{ x\in \widetilde{E}:\widetilde{r}%
<\left\| x\right\| _E<\widetilde{R}\right\}
\]
and its closure $\overline{D_{\widetilde{r},\widetilde{R}}}$.

Set $K:=\overline{B_{\widetilde{E}}\left( 0,\widetilde{r}\right) }\cup
\partial B_{\widetilde{E}}\left( 0,\widetilde{R}\right) $ and denote $%
K^\gamma :=\left\{ x\in E:\inf\limits_{z\in K}\left\| x-z\right\| _E<\gamma
\right\} $ for $\gamma >0$, where $B_{\widetilde{E}}\left( 0,\widetilde{r}%
\right) =\left\{ x\in \widetilde{E}:\left\| x\right\| _E<\widetilde{r}%
\right\} $, $\partial B_{\widetilde{E}}\left( 0,\widetilde{R}\right)
=\left\{ x\in \widetilde{E}:\left\| x\right\| _E=\widetilde{R}\right\} $.
Choose $\widetilde{r}$,$\delta >0$ sufficiently small such that
\[
\sup\limits_{x\in \overline{K^{2\delta }}}\Phi \left( x\right) <\frac{d^{*}}%
2,
\]
with $d^{*}:=\inf\limits_{x\in E,\text{ }\left\| x\right\| _E=r}\Phi \left(
x\right) $.

Let $\phi \in C^\infty \left( E,\left[ 0,1\right] \right) $ be the cut-off
function given by Proposition \ref{prop4.1}, satisfying $\phi \left(
x\right) =1$ for all $x\in K^\delta $ and $\phi \left( x\right) =0$ for all $%
x\in E\backslash K^{2\delta }$. Consider the follow $\eta :\Bbb{R}^{+}\times
E\rightarrow E$ defined by:
\[
\left\{
\begin{array}{cc}
\frac{d\eta \left( t,x\right) }{dt}=-\left( 1-\phi \left( \eta \left(
t,x\right) \right) \right) \Phi ^{\prime }\left( \eta \left( t,x\right)
\right) , & \left( t,x\right) \in R^{+}\times E, \\
\eta \left( 0,x\right) =x. &
\end{array}
\right.
\]
Define the minimax value
\[
d:=\inf\limits_{t\geq 0}\sup\limits_{u\in \overline{D_{\widetilde{r},%
\widetilde{R}}}}\Phi \left( \eta \left( t,x\right) \right) .
\]
Then $d$ is a critical value of $\Phi $ on $E$ with an upper bound
\begin{equation}
d\leq \widetilde{d}:=\sup\limits_{x\in \overline{D_{\widetilde{r},\widetilde{%
R}}}}\Phi \left( x\right) .  \label{eq92}
\end{equation}
Moreover, there exists a critical point $x^{*}\in E$ of $\Phi $ such that $%
\Phi \left( x^{*}\right) =d$ and its generalized Morse index satisfies $%
M\left( \Phi ,x^{*}\right) \geq l$.
\end{theorem}

Indeed, in Theorem \ref{thm4.5}, the $C^1$ smoothness of $\Phi $, together
with assumptions $\left\langle 1\right\rangle \left\langle 2\right\rangle $,
suffices to conclude that $d$ is a critical value of $\Phi $ on $E$, as the
proof of Theorem \ref{thm1.1} demonstrates. A concrete example is provided
by taking $E=X_n:=$span$\left\{ \varphi _1,\cdots ,\varphi _n\right\} $ and $%
\Phi =J$, where $\varphi _i$ are the eigenfunctions of the linear problem $%
\left( \ref{eq93}\right) $. The quantity $d$ defined in $\left( \ref{eq92}%
\right) $ depends on $l$ but is independent of $n$. The hypotheses of
Theorem \ref{thm1.1}, together with the additional assumption $J\in
C^3\left( X_n,\Bbb{R}\right) $(which holds for the mollified functionals $%
J_m $ on finite-dimensional subspaces of $X$), ensure that the full set of
assumptions in Theorem \ref{thm4.5} is satisfied. Therefore the assertion of
Theorem \ref{thm4.5} is valid in this setting. Moreover, since $%
\sup\limits_{x\in X_n}J\left( x\right) \rightarrow \infty $ as $n\rightarrow
\infty $, the critical point $x^{*}$ is a non-global saddle point of $J$ on $%
X_n$ for large $n$.

We need the following assumptions:

$\left\langle 2\right\rangle ^{*}$ $\Phi \left( x\right) \rightarrow -\infty
$ as $\left\| x\right\| _E\rightarrow \infty $;

$\left\langle 4\right\rangle $ $\Phi ^{\prime }\left( u\right) =u-F\left(
u\right) $, where $F:E\rightarrow E$ is a compact mapping.

Note that hypotheses $\left\langle 2\right\rangle $ and $\left\langle
3\right\rangle $ follow from $\left\langle 2\right\rangle ^{*}$, and the
(PS) condition follows from $\left\langle 2\right\rangle ^{*}$ together with
$\left\langle 4\right\rangle $. For the specific functional $J$ associated
with problem $\left( \ref{eq1}\right) $, the hypothesis $\left\langle
4\right\rangle $ (i.e., $J^{\prime }=I-F$ with $F$ compact) is satisfied due
to the subcritical growth condition $\left( f_1^{\prime }\right) $ and the
compact Sobolev embedding.

\begin{remark}
\label{rek4.6} Theorem \ref{thm4.5} provides an abstract setting for
generating critical points with high generalized Morse indices. However, the
condition $\dim E<\infty $ cannot be dropped. Indeed, for the case $\dim
E=\infty $, the four assumptions in Theorem \ref{thm4.5} --- namely, the $%
C^3 $ smoothness of $\Phi $, and the geometric conditions $\left\langle
1\right\rangle \left\langle 2\right\rangle ^{*}\left\langle 4\right\rangle $
--- cannot hold simultaneously. If they did, the theorem would produce a
sequence of critical points $\left\{ u_j\right\} $ of $\Phi $ with $M\left(
\Phi ,u_j\right) \geq j\rightarrow \infty $. However, from the mountain-pass
structure of $\Phi $(which follows from $\left\langle 1\right\rangle
\left\langle 2\right\rangle ^{*}$), if the sequence $\left\{ \left\|
u_j\right\| _E\right\} $ is unbounded, then we get a contradictory
inequality $d^{*}\leq \Phi \left( u_j\right) \rightarrow -\infty $.
Alternatively, if $\left\{ \left\| u_j\right\| _X\right\} $ is bounded, then
$u_j\rightarrow u_0^{*}\in K_\Phi :=\left\{ u\in E:\Phi ^{\prime }\left(
u\right) =0\right\} $ by (PS) condition. This also gives rise to a
contradiction $\infty >M\left( \Phi ,u_0^{*}\right) \geq M\left( \Phi
,u_j\right) \geq j\rightarrow \infty $.

Therefore, Theorem \ref{thm4.5} in its standard form for the case $\dim
E=\infty $ fails to hold. This highlights the necessity of the new minimax
framework developed in this paper, which does not require the original
functional $J$ to be globally $C^3$. Instead, we work with a sequence of
mollified $C^\infty $ functionals $J_m$ and exploit the smoothness within
finite-dimensional approximations, passing to the limit in the final step.
This allows us to produce critical points with arbitrarily high generalized
Morse indices under the assumption $\left( f_1^{\prime }\right) $, which is
indispensable for generalized Morse index analysis.
\end{remark}

\textit{Proof of Theorem 1.2}. In order to use Theorem \ref{thm1.1} to prove
Theorem \ref{thm1.2}, we need to verify that the nonlinear term satisfies
the assumptions of Theorem \ref{thm1.1}. By assumption $\left( f_3\right) $,
we have $f\left( x,t\right) =o\left( t\right) $ uniformly in $x$ as $%
t\rightarrow 0$. Since $2\leq q<p$, there exists a constant $C>0$ such that
\[
\left| f_t^{\prime }\left( x,t\right) \right| \leq C\left( 1+\left| t\right|
^{p-2}\right) .
\]
Take $\theta =q>2$ and $M>\left( \frac pq\right) ^{\frac 1{p-q}}$, then for $%
t\in \Bbb{R}$, $\left| t\right| \geq M$,
\[
0<\theta F\left( x,t\right) =\frac qp\left| t\right| ^p+\left| t\right|
^{q-1}t\leq \left| t\right| ^p+\left| t\right| ^{q-1}t=tf\left( x,t\right) ,
\]
where $F\left( x,t\right) =\int_0^tf\left( x,s\right) ds$. The assertion of
Theorem \ref{thm1.2} is validated.

Actually we can deal with a more general case by taking $f\left( x,t\right)
=a\left( x\right) \left| t\right| ^{p-2}t+b\left( x\right) \left| t\right|
^{q-1}$ instead of $f\left( x,t\right) =\left| t\right| ^{p-2}t+\left|
t\right| ^{q-1}$, and the conclusion is still established, where $a\left(
x\right) $, $b\left( x\right) \in C\left( \overline{\Omega }\right) $ and $%
a\left( x\right) >0$ on $\overline{\Omega }$.

\begin{remark}
\label{rek4.7} Let $0<\lambda _1<\lambda _2\leq \cdots \leq \lambda _j\leq
\cdots $ be the eigenvalues of the linear eigenvalue problem
\begin{equation}
\left\{
\begin{array}{cc}
-\Delta v=\lambda v, & x\in \Omega , \\
v=0, & x\in \partial \Omega ,
\end{array}
\right.  \label{eq93}
\end{equation}
and let $0<\mu _1<\mu _2\leq \cdots \leq \mu _j\leq \cdots $ be the
eigenvalues of the linear eigenvalue problem
\begin{equation}
\left\{
\begin{array}{cc}
-\Delta v+a\left( x\right) v=\mu v, & x\in \Omega , \\
v=0, & x\in \partial \Omega ,
\end{array}
\right.  \label{eq94}
\end{equation}
where $a\left( x\right) \in L^{\frac N2}\left( \Omega \right) $ and $N\geq 3$%
. The conclusions of Theorem \ref{thm1.1} remain valid for the following
generalized problems:

1. $-\Delta u=\lambda u+f\left( x,u\right) $ with $\lambda <\lambda _1$.

2. $-\Delta u+a\left( x\right) u=\mu u+f\left( x,u\right) $ with $\mu <\mu
_1 $.

In both cases, the nonlinearity $f$ is assumed to satisfy $\left(
f_1^{\prime }\right) $-$\left( f_3\right) $.
\end{remark}

State Key Laboratory of Mathematical Sciences, Academy of Mathematics and
Systems Science, Chinese Academy of Sciences, Beijing 100190, Peoples
Republic of China

Institute of Mathematics, Academy of Mathematics and Systems Science,
Chinese Academy of Sciences

School of Mathematical Sciences, University of Chinese Academy of Sciences,
Beijing 100049

Center for Excellence in Mathematical Sciences, Chinese Academy of Sciences

\textit{Email address}: \textrm{lichong@amss.ac.cn}

Institute of Mathematics, Academy of Mathematics and Systems Science,
Chinese Academy of Sciences

\textit{Email address}: \textrm{lisj@math.ac.cn}


\begin{thebibliography}{99}
\bibitem{AR}  Antonio Ambrosetti and Paul Rabinowitz. Dual variational
methods in critical point theory and applications. \textit{J.Funct.Anal.,}
14(4):349-381, 1973.

\bibitem{Ba}  Abbas Bahri. Topological results on a certain class of
functionals and applications. \textit{J.Funct.Anal.}, 41:379-427, 1981.

\bibitem{BCW}  Thomas Bartsch, Kungching Chang and Zhiqiang Wang. On the
Morse indices of sign changing solutions of nonlinear elliptic problems.
\textit{Math.Z.}, 233(4):655-677, 2000.

\bibitem{BW}  Thomas Bartsch and Tobias Weth. Three nodal solutions of
singularly perturbed elliptic equations on domains without topology. \textit{%
Ann.Inst.H.Poincar\'e-Anal.Non Lin\'eaire.,} 22(3):259-281, 2005.

\bibitem{BR}  Vieri Benci and Paul Rabinowitz. Critical point theorems for
indefinite functionals. \textit{Invent.Math.,} 52:241-273, 1979.

\bibitem{Br}  Haim Brezis. \textit{Functional Analysis}. Sobolev Spaces and
Partial Diffrential Equations.\emph{\ }Universitext, Springer, New York,
2011.

\bibitem{Cha1}  Kungching Chang. \textit{Infinite Dimensional Morse Theory
and Multiple Solution problems}. Birkhauser, Boston, 1993.

\bibitem{Cha2}  Kungching Chang. A variant mountain pass lemma. \textit{%
Scient.Sinica A.},\emph{\ }26:1241-1255, 1983.

\bibitem{EG}  Ivar Ekeland and Nassif Ghoussoub. $Z_2$-pequivariant
Ljusternik-Schnirelman theory for non-even functionals. \textit{%
Ann.Inst.H.Poincar\'e-Anal.Non Lin\'eaire., }15(3):341-370, 1998.

\bibitem{FG}  Guangcai Fang and Nassif Ghoussoub. Morse-type information on
Palais-Smale sequences obtained by min-max Principles. \textit{%
Comm.Pure.Appl.Math.}, 47(12):1595-1653, 1994.

\bibitem{Ghou}  Nassif Ghoussoub. \textit{Duality and Perturbation Methods
in Critical Point Theory.} Cambridge Tracts in Mathematics. Cambridge
University Press, 1993.

\bibitem{Lang}  Serge Lang. \textit{Differential Manifolds}. Second edition,
Springer-Verlag, New-York, 1985.

\bibitem{LiLi}  Chong Li and Shujie Li. Gaps of consecutive eigenvalues of
Laplace operator and the existence of multiple solutions for superliner
elliptic problem. \textit{J.Funct.Anal.}, 271(1):245-263, 2016.\emph{\ }

\bibitem{LiLiu}  Chong Li and Yanyan Liu. Multiple solutions for a class of
semilinear elliptic problems via Nehari-type linking theorem. \textit{%
Calc.Var.Partial Differential Equations.}, 56(2):1-14, 2017.

\bibitem{LLZ}  Shujie Li and Zhaoli Liu. Perturbations from Symmetric
Elliptic Boundary Value Problems. \textit{J.Differential Equations.},\emph{\
}185(1):271--280, 2002.

\bibitem{LW}  Shujie Li and Michel Willem. Applications of local linking to
critical point theory. \textit{J.Math.Anal.Appl.}, 189(1):6-32, 1995.

\bibitem{Ne}  Zeev Nehari. Characteristic values associated with a class of
nonlinear second-order differential equations. \textit{Acta Math.},
105:141-175, 1961.

\bibitem{Ni}  Weiming Ni. Some minimax principle and their applications in
nonlinear elliptic equations. \textit{J.d'Analyse Math.}, 37(1):248-275,
1980.\emph{\ }

\bibitem{RSW}  Paul H. Rabinowitz, Jiabao Su and Zhiqiang Wang. Multiple
solutions of superlinear elliptic equations. \textit{Rend. Lincei Mat.Appl.}%
, 18(1):97-108, 2007.

\bibitem{Ra1}  Paul H. Rabinowitz. Some critical point theorems and
applications to semilinear elliptic partial differential equations. \textit{%
Ann.Sc.Norm.Sup.Pisa,Ser.}, 4(1):215-223, 1978.

\bibitem{Ra2}  Paul H. Rabinowitz. Some aspects of Nonlinear eigenvalue
problems. \textit{Rocky Mountain Journal of Mathematics.}, 3(2):161-201,
1973.

\bibitem{Ra3}  Paul H. Rabinowitz. Multiple critical points of perturbed
symmetric functionals. \textit{Trans.Amer.Math.Soc.}, 272(2):753-769, 1982.

\bibitem{Ra4}  Paul H. Rabinowitz. \textit{Minimax metheds in critical point
theory with applications to differential equations}. CBMS, Regional
Conference Series in Mathematics, Number 65, 1986.

\bibitem{SST}  Mingzheng Sun, Jiabao Su and Rushun Tian. Five nontrivial
solutions of superlinear elliptic problem. \textit{J.Funct.Anal.},
286(8):110266, 2024.

\bibitem{Stru1}  Michael Struwe. Infinitely many critical points for
functionals which are not even and applications to superlinear boundary
value problems. \textit{Manusc.Math.}, 32:335-364, 1980.

\bibitem{Stru2}  Michael Struwe. Multiple solutions of anticoercive boundary
value problems for a class ordinary differential equations of second order.
\textit{J.Differential Equations.}, 37(3):285-295, 1980.

\bibitem{Stru3}  Michael Struwe. \textit{Variational methods and their
applications to non-linear partial differential equations and Hamiltonian
systems}. Springer-Verlag, 1996.

\bibitem{Szu}  Andrzej Szulkin. Ljusternik-Schnirelmann theory on $C^1$%
-manifold. \textit{Ann.Inst.H.Poincar\'e-Anal.Non Lin\'eaire.},
5(2):119-139, 1988.

\bibitem{Turner}  Robert E.L. Turner. Superlinear Sturm-Liouville problems.
\textit{J.Differential Equations.}, 13(1):157-171, 1973.

\bibitem{W}  Zhiqiang Wang. On a superlinear elliptic equation. \textit{%
Ann.Inst.H.Poincar\'e-Anal.Non Lin\'eaire.}, 8(1):43-58, 1991.

\bibitem{Willem}  Michel Willem. \textit{Minimax theorems}. Progress in
Nonlinear Differential Equations and Their Applications, 24. Birkh\"auser,
Boston, 1996.
\end{thebibliography}
\end{document}